\documentclass{gtpart}

\gtart

\usepackage{amscd,latexsym,graphicx}
\usepackage[all]{xy}
\usepackage{enumerate}

\title{$K$-duality for stratified pseudomanifolds}

\author{Claire Debord}
\givenname{Claire}
\surname{Debord}
\address{Laboratoire de Math{\'e}matiques, Universit{\'e} Blaise Pascal,
  Complexe universitaire dse C{\'e}zeaux, 24 Av. des Landais, 63177
  Aubi{\`e}re cedex, France}
\email{debord@math.univ-bpclermont.fr}
\urladdr{http://math.univ-bpclermont.fr/~debord/}

\author{Jean-Marie Lescure}
\givenname{Jean-Marie}
\surname{Lescure}
\address{Laboratoire de Math{\'e}matiques, Universit{\'e} Blaise Pascal,
  Complexe universitaire dse C{\'e}zeaux, 24 Av. des Landais, 63177
  Aubi{\`e}re cedex, France}
\email{lescure@math.univ-bpclermont.fr}
\urladdr{http://math.univ-bpclermont.fr/~lescure/}

\keyword{Singular manifolds}
\keyword{Smooth groupoids}               
\keyword{Kasparov bivariant $K$-theory}
\keyword{Poincar{\'e} duality}
\subject{primary}{msc2000}{58B34, 46L80, 19K35, 58H05,  57N80}
\subject{secondary}{msc2000}{19K33, 19K56, 58A35 }

\volumenumber{}
\issuenumber{}
\publicationyear{}
\papernumber{}
\startpage{}
\endpage{}
\doi{}
\MR{}
\Zbl{}
\received{}
\revised{}
\accepted{}
\published{}
\publishedonline{}
\proposed{}
\seconded{}
\corresponding{}
\editor{}
\version{}

\newtheorem{theorem}{Theorem}
\newtheorem{proposition}{Proposition}
\newtheorem{corollary}{Corollary}
\newtheorem{lemma}{Lemma}

\theoremstyle{definition}
\newtheorem{definition}{Definition}

\theoremstyle{definition}
\newtheorem{examples}{Examples}

\theoremstyle{definition}
\newtheorem{example}{Example}

\theoremstyle{definition}
\newtheorem{remark}{Remark}

\theoremstyle{definition}

\theoremstyle{definition}

\newcommand{\cA}{\mathcal{A}}

\newcommand{\cG}{\mathcal{G}}

\newcommand{\cI}{\mathcal{I}}

\newcommand{\cK}{\mathcal{K}}

\newcommand{\cN}{\mathcal{N}}

\newcommand{\cS}{\mathcal{S}}

\newcommand{\fS}{ \mathsf{S}}
\newcommand{\fR}{ \mathsf{R}}

\newcommand{\NN}{\mathbb N}
\newcommand{\ZZ}{\mathbb Z}

\newcommand{\RR}{\mathbb R}
\newcommand{\CC}{\mathbb C}

\newcommand{\im}{\mathop{\mathrm{Im}}\nolimits}

\newcommand{\id}{\mathop{\mathrm{Id}}\nolimits}

\newcommand{\ot}{\otimes}

\newcommand{\pb}[1]{{}^*#1^*}

\newcommand\bfx{\mathbf{x}}
\newcommand\bfv{\mathbf{v}}
\newcommand\bfw{\mathbf{w}}

\newcommand{\wtU}{\widetilde{U}}
\newcommand{\wtphi}{\widetilde{\phi}}
\newcommand{\wtpsi}{\widetilde{\psi}}

\newcommand{\bU}{\overline{U}}
\newcommand{\bV}{\overline{V}}

\newcommand\smx{X^{\circ}}
\newcommand\smX{X^{\circ}}
\newcommand\smA{A^{\circ}}

\newcommand\smF{F^{\circ}}
\newcommand\smO{O^{\circ}}
\newcommand\sml{L^{\circ}}

\newcommand\smN{\cN^{\circ}}
\newcommand\pa{\partial}

\newcommand{\interior}[1]{\overset{\circ}{#1}}

\numberwithin{equation}{section}

\begin{document}

\begin{abstract} This paper continues the project started in \cite{DL2} where Poincar\'e
duality in $K$-theory was studied for singular manifolds with
isolated conical singularities. Here,  we extend the study and the
results to general stratified pseudomanifolds.  We review the
axiomatic definition of a smooth stratification $\fS$ of a topological
space $X$ and we define a groupoid $T^{\fS}X$,  called the $\fS$-tangent space. 
This groupoid is made of different pieces encoding the tangent
spaces of strata,  and these pieces are glued into the smooth
noncommutative groupoid $T^{\fS}X$ using the familiar procedure
introduced by A. Connes for the tangent 
groupoid of a manifold. The main result is that $C^{*}(T^{\fS}X)$ is
Poincar\'e dual to $C(X)$,  in other words, the $\fS$-tangent space
plays the role in $K$-theory  of a tangent space for $X$.   
\end{abstract}

\maketitle


\section*{Introduction}
This paper takes place in a longstanding project aiming to study index
theory and related questions on stratified pseudomanifolds using tools
and concepts from noncommutative geometry. 

The key observation at the beginning of this project is that in its
$K$-theoritic form,  the Atiyah-Singer index
theorem \cite{AS1} involves ingredients that should survive to the
singularities allowed in a stratified pseudomanifold. This is possible, from
our opinion, as soon as one accepts
reasonable generalizations and new presentation of certain classical
objects on smooth manifolds, making sense on stratified pseudomanifolds. 

The first instance of these classical objects that need to be adapted
to singularities is the notion of {\sl tangent space}. Since index
maps in \cite{AS1} are defined on the $K$-theory of the tangent spaces
of smooth manifolds, one must have a similar space adapted to
stratified pseudomanifolds. Moreover, such a space should satisfy
natural requirements. It should coincide with the usual notion on the
regular part of the pseudomanifold and incorporate in some way copies
of usual tangent spaces of strata, while keeping enough smoothness to
allow interesting computations. Moreover, it should be Poincar{\'e} dual
in $K$-theory (shortly, $K$-dual) to the pseudomanifold itself. This
$K$-theoritic property involves bivariant $K$-theory and was proved
between smooth manifolds and their tangent spaces by G. Kasparov
\cite{Ka1} and A. Connes-G. Skandalis \cite{AC-GS1984}. 

In \cite{DL2}, we
introduced a candidate to be the tangent space of a pseudomanifold
with isolated conical singularities. It appeared to be a smooth
groupoid, leading to a noncommutative $C^*$-algebra, and we proved
that it fulfills the expected 
$K$-duality. 

In \cite{JM_L2006}, the second author interpreted the duality proved in
\cite{DL2} as a principal symbol map,  thus recovering the  classical
picture of Poincar\'e duality in $K$-theory for smooth
manifolds. This interpretation used a notion of noncommutative
elliptic symbols,  which appeared to be the cycles of the 
$K$-theory of the noncommutative tangent space. 

In \cite{DLN2006}, the noncommutative tangent space together with other
deformation groupoids was used to construct analytical and topological
index maps, and their equality was proved. As expected, these index
maps are straight generalizations of those of \cite{AS1} for
manifolds. 

The present paper is devoted to the construction of the noncommutative
tangent space for a general stratified pseudomanifold and the proof of
the $K$-duality. It is thus a sequel of \cite{DL2}, but can be read
independently. At first glance, one should have expected that the
technics of \cite{DL2} iterate easily to give the general result. In
fact, although the definition of the groupoid giving the noncommutative
tangent space itself is natural and intuitive in the general case, its
smoothness is quite intricate and brings issues that did not exist in
the conical case. We have given here a detailed treatment of this
point, since we believe that this material will be usefull in further
studies about the geometry of stratified spaces. Another difference
with \cite{DL2} is that we have given up the explicit construction of a
dual Dirac element. Instead, we use an easily defined Dirac element
and then prove the Poincar{\'e} duality by an induction, based on
an operation called unfolding which consists in removing the minimal strata in a
pseudomanifold and then ``doubling'' it to get a new pseudomanifold,
less singular. The difficulty in this approach is moved to the proof
of the commutativity of certain diagrams in $K$-theory, necessary to
apply the five lemma and  to continue the induction. 

The interpretation of this $K$-duality in terms of noncommutative
symbols and pseudodifferential operators, as well as the construction
of index maps together with the statement of an index theorem, is postponed
to forthcoming papers.

This approach of index theory on singular spaces in the framework of
noncommutative geometry takes place in a long history of past and
present resarch works. But the specific issues about Poincar{\'e}
duality, bivariant $K$-theory, topological index maps and statement of
 Atiyah-Singer like theorems are quite recent and attract an increasing
interest \cite{MN2005,MR2006,EM2005,EM2007,Sav1,NSS2}. 

\subsection*{Acknowledgements} We would like to thank the referee for
his useful comments and remarks. In particular he suggested to us
Lemmas \ref{referee-1} and \ref{referee-2} which enable us to shorten and clarify
significantly the proof of Theorem \ref{main-theorem}.


\section{Basic definitions}

\subsection{Around Lie groupoids} \label{groupoids}

We refer to \cite{Re,CW,Mc,DL2007} for the classical definitions and
constructions related to groupoids, their Lie algebroids and 
$C^*$-algebras of groupoids. In this section, we fix the notations and recall the
less classical definitions and results needed in the sequel. Some
material presented here is already in \cite{DL2,DLN2006}.

\subsubsection{Pull back groupoids}
Let $G\rightrightarrows M$ be a locally compact Hausdorff groupoid with source $s$ and range
$r$. If $f:N\rightarrow M$ is a surjective map, the {\it pull
back} groupoid $\pb{f}(G)\rightrightarrows N$ of $G$ by $f$ is by
definition the set
\begin{equation*}
\pb{f}(G):=\{(x,\gamma,y)\in N\times G\times N \ \vert \
r(\gamma)=f(x),\ s(\gamma)=f(y)\}
\end{equation*}
with the structural morphisms given by
\begin{enumerate}
\item the unit map $x \mapsto (x,f(x),x)$,

\item the source map $(x,\gamma,y)\mapsto y$ and range map
$(x,\gamma,y) \mapsto x$,

\item the product $(x,\gamma,y)(y,\eta,z)=(x,\gamma \eta ,z)$ and
inverse $(x,\gamma,y)^{-1}=(y,\gamma^{-1},x)$.
\end{enumerate}

\smallskip \noindent The results of  \cite{MRW} apply to show that the groupoids $G$
and $\pb{f}(G)$ are Morita equivalent when $f$ is surjective and open.

\smallskip \noindent Let us assume for the rest of this subsection that $G$ is a smooth
groupoid and that $f$ is a surjective submersion, then $\pb{f}(G)$
is also a Lie groupoid. Let $(\cA(G), q, [\ ,\ ])$ be the Lie
algebroid of $G$. Recall
that $q: \cA(G) \to TM$ is the anchor map. Let
$(\cA(\pb{f}(G)),p,[\ ,\ ])$ be the Lie algebroid of $\pb{f}(G)$
and $Tf : TN \to TM$ be the differential of $f$. Then there exists an isomorphism 
\begin{equation*}
    \cA(\pb{f}(G)) \simeq \{(V, U)\in TN \times \cA(G) \
    \vert \ Tf(V)=q(U) \in TM \}
\end{equation*}
under which the anchor map $p:\cA(\pb{f}(G)) \rightarrow TN$
identifies with the projection $TN \times \cA(G) \to TN$. (In
particular, if $(V, U) \in \cA(\pb{f}(G))$ with $V \in T_xN$ and
$U \in \cA_y(G)$, then $y = f(x)$.)

\subsubsection{Subalgebras and exact sequences of groupoid 
  $C^*$-algebras} 

To any smooth groupoid $G$ are associated two
$C^*$-algebras 
corresponding to
two different completions of the
involutive convolution algebra $C_c^\infty 
(G)$, namely  the reduced and maximal $C^*$-algebras \cite{Co3,Co0,Re}. We will denote
respectively these
$C^*$-algebras by $C^*_{r}(G)$ and $C^*(G)$. Recall that the
identity on $C_c^\infty 
(G)$ induces a surjective morphism from $C^*(G)$ onto
$C^*_{r}(G)$ which is an isomorphism if the
groupoid $G$ is amenable. Moreover in this case the $C^*$ algebra of
$G$ is nuclear \cite{AR2000} . 

\smallskip \noindent We will
use the following usual notations:\\ 
Let $G \underset{r}{\overset{s}{\rightrightarrows}} G^{(0)}$ be a smooth 
groupoid with source $s$ and range $r$. 
If $U$ is any subset of $G^{(0)}$, we let: 
$$G_U:=s^{-1}(U)\ , \ G^U:=r^{-1}(U) \mbox{ and } 
G_U^U=G\vert_U:=G_U\cap G^U \ .$$ 

\smallskip \noindent To  an open subset  $O$  of $G^{(0)}$ corresponds an inclusion 
$i_O$ of $C_c^\infty(G\vert_O)$ into $C_c^\infty(G)$ which induces an 
injective morphism, again denoted by $i_O$, from  $C^*(G\vert_O)$ into 
$C^*(G)$. \\ When $O$ is saturated, $C^*(G\vert_O)$ 
is an ideal of $C^*(G)$. In this case, $F:=G^{(0)}\setminus O$ 
is a saturated closed subset of $G^{(0)}$ and the restriction of 
functions induces a  surjective morphism $r_F$ from $C^*(G)$ to 
$C^*(G\vert_F)$. Moreover, according to \cite{HS1}, the following sequence 
of $C^*$-algebras is exact:
$$\begin{CD} 0 @>>> C^*(G\vert_O) @>{i_O}>> C^*(G) 
@>{r_F}>>  C^*(G\vert_F) @>>> 0 \end{CD} \ .$$

\subsubsection{$KK$-elements associated to deformation groupoids}  
 
A smooth groupoid $G$ is called a {\it deformation groupoid} 
if: 
$$G= G_1 \times \{0\} \cup G_2\times ]0,1] 
\rightrightarrows G^{(0)}=M\times [0,1],$$ 
where $G_1$ and $G_2$ are smooth groupoids with unit space $M$. That 
is, $G$ is obtained by gluing 
$G_2\times ]0,1]\rightrightarrows M\times ]0,1]$, 
which is the cartesian product of the groupoid
  $G_2\rightrightarrows M$ with the space $]0,1]$, with the groupoid 
$G_1\times\{0\}\rightrightarrows M\times\{0\}$. 
 
\smallskip \noindent In this situation one can consider the saturated 
open subset $M\times ]0,1]$ of $G^{(0)}$. Using the isomorphisms 
$C^*(G\vert_{M\times ]0,1]}) \simeq C^*(G_2)\ot C_0(]0,1])$ and 
$C^*(G\vert_{M\times\{0\}})\simeq C^*(G_1)$,  we obtain the following 
exact sequence of $C^*$-algebras: 
$$\begin{CD} 0 @>>> C^*(G_2)\ot C_0(]0,1]) @>{i_{M\times]0,1]}}>> C^*(G) 
@>{ev_0}>>  C^*(G_1)  @>>> 0 \end{CD} 
$$ 
where $i_{M\times]0,1]}$ is the inclusion map and $ev_0$ is the {\it evaluation map} at 
$0$, that is $ev_0$ is the map coming from the 
restriction of functions to $G\vert_{M\times\{0\}}$.

\noindent We assume now that $C^*(G_1)$ is nuclear. Since 
the $C^*$-algebra $C^*(G_2)\ot C_0(]0,1])$ is contractible, 
 the long exact sequence in $KK$-theory shows that the group homomorphism 
$(ev_0)_*=\cdot {\ot}[ev_0]:KK(A,C^*(G)) \rightarrow KK(A, C^*(G_1))$ 
is an isomorphism for each $C^*$-algebra $A$ \cite{Ka1,JC-GS1986}. 
 
\noindent In particular with $A=C^*(G_1)$ and $A=C^*(G)$
we get that $[ev_0]$ is invertible in 
$KK$-theory: there is an element $[ev_0]^{-1}$ in 
$KK(C^*(G_1), C^*(G))$ such that $[ev_0]^{-1} {\ot}
[ev_0]=1_{C^*(G_1)}$ 
and $[ev_0] {\ot} [ev_0]^{-1}=1_{C^*(G)}$. 
 
\smallskip \noindent Let $ev_1:C^*(G) \rightarrow C^*(G_2)$ be the 
evaluation map at $1$ and $[ev_1]$ the corresponding element of 
$KK(C^*(G),C^*(G_2))$. 

\smallskip \noindent The {\it $KK$-element associated to the 
  deformation groupoid} $G$ is defined by: 
$$\delta=[ev_0]^{-1} {\ot} [ev_1]\in  KK(C^*(G_1),C^*(G_2)) \ . $$ 

\noindent One can find examples of such elements related to index theory in
\cite{Co0,HS1,DL2,DLN2006,DL2007}.

\subsection{Generalities about $K$-duality}

We give in this paragraph some general facts about Poincar{\'e} duality in
bivariant K-theory. Most of them are well known and proofs are only added
when no self contained proof could be found in
the literature. All $C^*$-algebras are assumed to be separable and
$\sigma$-unital. 

\smallskip \noindent Let us first recall what means the Poincar{\'e}
duality in $K$-theory \cite{Ka2,AC-GS1984,Co0}: 
\begin{definition}\label{def-Kduality}
  Let $A,B$ be two $C^*$-algebras. One says that $A$ and $B$
  are Poincar{\'e} dual, or shortly $K$-dual, when there exists 
  $\alpha\in K^0(A\otimes B)=KK(A\otimes B, \CC)$ and 
  $\beta\in KK(\CC,A\otimes B)\simeq K_0(A\otimes B)$ 
  such that 
   $$
     \beta\underset{B}{\otimes}\alpha =1\in KK(A,A) \hbox{ and }  
     \beta\underset{A}{\otimes}\alpha =1\in KK(B,B)
   $$  
Such elements are then called  Dirac and dual-Dirac elements. 
\end{definition}

\noindent It follows that for $A,B$ two $K$-dual $C^*$-algebras and for any
$C^*$-algebras $C,D$, the following
isomorphisms hold: 
 $$ 
   \beta\underset{B}{\otimes}\cdot: KK(B\otimes C, D)
    \longrightarrow KK(C,A\otimes D);
 $$
 $$
   \beta \underset{A}{\otimes}\cdot: KK(A\otimes C, D)
    \longrightarrow KK(C,B\otimes D);
 $$
with inverses given respectively by $\cdot\underset{A}{\otimes}\alpha$
and $\cdot\underset{B}{\otimes}\alpha$. 

\begin{example}A basic example is $A=C(V)$ and $B=C_0(T^*V)$ where $V$ is a
 closed smooth manifold (\cite{Ka2,AC-GS1984}, see also \cite{DL2} for a
 description of the Dirac element in terms of groupoids). This duality
 allows to recover that the usual quantification and principal
 symbol maps are mutually inverse isomorphisms in $K$-theory: 
$$\Delta_V =(\cdot \underset{C_0(T^*V)}{\otimes} \alpha):
{K_0(C_0(T^*V))} \overset{\simeq}{\longrightarrow} {K^0(C(V))} $$
$$\Sigma_V =(\beta \underset{C(V)}{\otimes} \cdot):
{K^0(C(V))} \overset{\simeq}{\longrightarrow} {K_0(C_0(T^*V))} $$
\end{example}

\noindent We observe that: 
\begin{lemma}\label{Kdual}
  Let $A,B$ be two $C^*$-algebras. Assume that there exists 
  $\alpha\in KK(A\otimes B,\CC)$ and 
  $\beta,\beta'\in KK(\CC,A\otimes B)$ satisfying 
  $$
    \beta\underset{B}{\otimes}\alpha = 1 \in KK(A,A) \hbox{ and } 
    \beta'\underset{A}{\otimes}\alpha = 1 \in KK(B,B) 
  $$
Then $\beta=\beta'$ so $A,B$ are $K$-dual. 
\end{lemma}
\begin{proof}
A simple calculation shows that for all $x\in KK(C,A\otimes D)$ we
have: $$ \beta\underset{B}{\otimes}(x\underset{A}{\otimes}\alpha) =x
\underset{A\otimes B}{\otimes}(\beta \underset{B}{\otimes} \alpha ) \ .$$
Applying this to $C=\CC$, $D=A$ and $x=\beta'$ we get: 
  $$
   \beta'=
   \beta\underset{B}{\otimes}(\beta'\underset{A}{\otimes}\alpha)
   =\beta\underset{B}{\otimes}1 = \beta 
  $$
\end{proof}
\begin{corollary}\label{K-duality-by-Dirac}
1) Given two $K$-dual $C^*$-algebras and
a Dirac element $\alpha$, the dual-Dirac element $\beta$ satisfying
the definition \ref{def-Kduality} is unique. \\
2) If there exists 
  $\alpha\in KK(A\otimes B, \CC)$ such that  
  $$ 
   \cdot \underset{B}{\otimes}\alpha: KK(\CC, A\otimes B)\longrightarrow
   KK(A,A) \hbox{ and } \cdot \underset{A}{\otimes}\alpha: KK(\CC, A\otimes B)\longrightarrow
   KK(B,B)
  $$
  are onto, then $A,B$ are $K$-dual and $\alpha$ is a Dirac element.
\end{corollary}
The two lemmas below have been comunicated to us by the referee. 
\begin{lemma}\label{referee-1}
 Let $J_1$ and $J_2$ be two closed two sided ideals in a nuclear
 $C^*$-algebra $A$ such that $J_1\cap J_2=\{0\}$  
and set $B=A/(J_1+J_2)$. Denote by $\partial_k\in KK_1(B,J_k)$,
$k=1,2$, the $KK$-elements associated respectively with the exact
sequences $0\longrightarrow J_1 \longrightarrow A/J_2\longrightarrow
B\longrightarrow 0$ and  
$0\longrightarrow J_2 \longrightarrow A/J_1\longrightarrow
B\longrightarrow 0$. Let also $i_k: J_k\longrightarrow A$ denote the
inclusions. Then the following equality holds:  
    $$(i_1)_*(\partial_1)+(i_2)_*(\partial_2)=0.$$
\end{lemma}
\begin{proof}
 Let $\partial\in KK_1(B,J_1+J_2)$ denote the $KK$-element associated with the exact sequence 
$0\longrightarrow J_1+J_2 \longrightarrow A\longrightarrow
B\longrightarrow 0$. Denote by $j_k:J_k\longrightarrow J_1+J_2$  
and $i:J_1+J_2\longrightarrow A$ the inclusions and by $p_k:
J_1+J_2\longrightarrow J_k$ the projections, $k=1,2$. Since the diagrams $(k=1,2)$
\begin{equation}
  \xymatrix{
   0 \ar[r] & J_1+J_2 \ar[d]_{p_k} \ar[r] & A \ar[d]\ar[r] & B\ar[r]\ar[d]_{=} & 0 \\
   0 \ar[r] & J_k\ar[r] & A/J_{3-k}\ar[r] & B\ar[r] & 0 } 
\end{equation}
commute, it follows that $(p_k)_*(\partial)=\partial_k$. \\ Moreover, 
$(p_1)_*\times(p_2)_*: KK_1(B,J_1+J_2)\longrightarrow
KK_1(B,J_1)\times KK_1(B,J_2)$ is an isomorphism whose inverse is $(j_1)_*+(j_2)_*$. It follows that 
$\partial = (j_1)_*(\partial_1)+(j_2)_*(\partial_2)$. Moreover the
six-term exact sequence associated to $0\longrightarrow J_1+J_2 \longrightarrow A\longrightarrow
B\longrightarrow 0$ leads to $i_*(\partial)=0$. The result follows now
from the equalities $i_k=i\circ j_k$, $k=1,2$.
\end{proof}
\begin{lemma}\label{referee-2}
 Let $X$ be a compact space and $A$ be a nuclear $C(X)$-algebra. Let $U_1$ and $U_2$ be disjoint open subsets of $X$. Set 
$X_1=X\setminus U_2$ and $J_k=C_0(U_k)A$, $k=1,2$. Let $\Psi : C(X)\otimes A\longrightarrow A$ be the homomorphism defined by 
$\Psi(f\otimes a)=fa$ and let $\varphi : C(X_1)\otimes J_1\longrightarrow A$, 
$\psi : C_0(U_2)\otimes A/J_1\longrightarrow A$
be the homomorphisms induced by $\Psi$.\\
Denote by $\overline{\partial_1}\in KK_1(A/J_1,J_1)$ and $\widetilde{\partial_2}\in KK_1(C(X_1),C_0(U_2))$ the $KK$-elements associated respectively with the exact sequences 
$0\longrightarrow J_1 \longrightarrow A\longrightarrow A/J_1\longrightarrow 0$ and
$0\longrightarrow C_0(U_2) \longrightarrow C(X)\longrightarrow
C(X_1)\longrightarrow 0$. Then the following equality holds:  
  $$(\varphi)_*(\overline{\partial_1}\otimes 1_{C(X_1)})+(\psi)_*(1_{A/J_1}\otimes\widetilde{\partial_2})=0.$$
\end{lemma}
\begin{proof}
We use the notation of Lemma \ref{referee-1}. We have commuting diagrams
\begin{equation}
  \xymatrix{
   0 \ar[r] & J_1\otimes C(X_1) \ar[d]_{\varphi_1}  \ar[r] & A \otimes
   C(X_1)\ar[d]\ar[r] & A/J_1\otimes C(X_1)\ar[r] \ar[d]_{\chi} & 0 \\
   0 \ar[r] & J_1\ar[r] & A/J_{2}\ar[r] & B\ar[r] & 0 } 
\end{equation}
and
\begin{equation}
  \xymatrix{
   0 \ar[r] & A/J_1\otimes C_0(U_2) \ar[d]_{\psi_2} \ar[r] & A/J_1\otimes C(X) \ar[d]\ar[r] & A/J_1\otimes C(X_1)\ar[r] \ar[d]_{\chi}& 0 \\
   0 \ar[r] & J_2\ar[r] & A/J_{1}\ar[r] & B\ar[r] & 0 } 
\end{equation}
where vertical arrows are induced by $\Psi$. It follows that
$(\varphi_1)_*(\overline{\partial_1}\otimes
1_{C(X_1)})=\chi^*(\partial_1)$ and $(\psi_2)_*(1_{A/J_1}\otimes
\widetilde{\partial_2})=\chi^*(\partial_2)$. We then use the equalities
$\varphi=i_1\circ \varphi_1$ and $\psi=i_2\circ \psi_2$ and 
apply Lemma \ref{referee-1} to conclude. 
\end{proof}
It yields the following, with the notation of Lemma \ref{referee-2} : 
\begin{lemma}\label{2-ou-of-3}
Let $\delta$ be in $K^0(A)$ and set $D=\Psi^*(\delta)$, 
$D_1=\varphi^*(\delta)$, $D_2=\psi^*(\delta)$. Then for any $C^*$-algebras $C$ and $D$, the two following 
long diagrams commute:
{\scriptsize \begin{equation} \label{2sur3-1} 
\xymatrix{
  \cdots KK_i(C,D\! \otimes \! C_0(U_2))
  \ar[r]\ar[d]^{\!\underset{C_0(U_2)}{\otimes}c_{i} D_2}
  & KK_i(C,D\! \otimes \!C(X))\ar[r]\ar[d]^{\underset{C(X)}{\otimes}c_{i}D} & 
 KK_i(C,D\! \otimes\! C(X_1))\ar[r]\ar[d]^{\underset{C(X_1)}{\otimes}c_{i}D_{1}}& KK_{i+1}(C,D\!\otimes\!
 C_0(U_2))\ar[d]^{\!\underset{C_0(U_2)}{\otimes}c_{i+1}D_{2}}\cdots\\
 \cdots KK_i(C\!\otimes\! A/J_1, D)\ar[r] & KK_i(C\!\otimes\!
  A,D)\ar[r]& KK_i(C\!\otimes\! J_1,D)\ar[r]&
  KK_{i+1}(C\!\otimes\! A/J_1,D) \cdots}
\end{equation}}
{\scriptsize \begin{equation} \label{2sur3-2} 
\xymatrix{
 \cdots KK_i(C,D\!\otimes\! J_1) \ar[r]\ar[d]^{\underset{J_1}{\otimes}c_{i}D_1} & 
  KK_i(C,D\!\otimes\! A)\ar[r]\ar[d]^{\underset{A}{\otimes}c_{i}D} & 
 KK_i(C,D\!\otimes\! A/J_1)\ar[r]\ar[d]^{\underset{A/J_1}{\otimes}c_{i}D_{2}}& KK_{i+1}(C,D\!\otimes\!
 J_1)\ar[d]^{\underset{J_1}{\otimes}c_{i+1}D_{1}}\cdots\\
 \cdots KK_i(C\!\otimes\! C(X_1),D) \ar[r] & KK_i(C\!\otimes\! C(X),D)\ar[r]& KK_i(C\!\otimes \! C_0(U_2),D)\ar[r]&
  KK_{i+1}(C\!\otimes \!C(X_1),D)\cdots}
\end{equation}}
where the $c_i$ belong to $\{-1,1\}$ and are chosen such that
$c_i=(-1)^{i+1}c_{i+1}$.\\ In particular, if
two of three elements $D_1,D_2,D$ are Dirac elements, so is the
third one. 
\end{lemma}
\begin{proof}
Observe first that Lemma \ref{referee-2} reads: 
$\overline{\partial_1}\underset{J_1}{\otimes}[\varphi]=
-\widetilde{\partial_2}\underset{C_0(U_2)}{\otimes}[\psi]$, 
which gives:
 $$
 \overline{\partial_1}\underset{J_1}{\otimes} D_1
=\overline{\partial_1}\underset{J_1}{\otimes}([\varphi]\otimes \delta) 
=(\overline{\partial_1}\underset{J_1}{\otimes}[\varphi])\otimes \delta
=(-\widetilde{\partial_2}\underset{C_0(U_2)}{\otimes}[\psi])\otimes \delta
=-\widetilde{\partial_2}\underset{C_0(U_2)}{\otimes}D_2
 $$
Now, using the skew-commutativity of the product $\underset{\CC}{\otimes}$, we
have for any $x\in KK_i(C,D\! \otimes\! C(X_1))$: 
\begin{eqnarray*}
 \overline{\partial_1}\underset{J_1}{\otimes}( x\underset{C(X_1)}{\otimes} D_1)
 &=& (\overline{\partial_1}\underset{\CC}{\otimes} x)\underset{J_1\otimes C(X_1)}{\otimes} D_1\\
 &=& (-1)^{i}(x \underset{\CC}{\otimes}\overline{\partial_1})\underset{J_1\otimes
     C(X_1)}{\otimes} D_1 \\
 &=& (-1)^{i} x \underset{C(X_1)}{\otimes} (\overline{\partial_1}\underset{J_1}{\otimes} D_1)\\
 &=& (-1)^{i} x \underset{C(X_1)}{\otimes} (-\widetilde{\partial_2}\underset{C_0(U_2)}{\otimes}D_2)\\
 &=& (-1)^{i+1}  (x \underset{C(X_1)}{\otimes}\widetilde{\partial_2})\underset{C_0(U_2)}{\otimes}D_2
\end{eqnarray*}
This yields,  thanks to the choice of the sign $c_{i}$, the
commutativity for the squares involving boundary homomorphisms in
Diagram (\ref{2sur3-1}).  
The other squares in Diagram (\ref{2sur3-1}) commute by 
definition of $D_{1}, D_{2}, D$ and by functoriality of
$KK$-theory. The commutativity of Diagram (\ref{2sur3-2}) is proved
by the same arguments.  The last assertion is then a 
consequence of Corollary \ref{K-duality-by-Dirac} and the five lemma. 
\end{proof}


\section{Stratified pseudomanifolds}

We are interested in studying stratified pseudomanifolds 
\cite{W1965,Mat,GoMa}. We will use the notations and equivalent descriptions
given by A. Verona in \cite{Ver} or used by J.P. Brasselet, G.
Hector and M. Saralegi in \cite{BHS}. The reader should
also look at \cite{HW} for a hepfull survey of the subject.

\subsection{Definitions} 

Let $X$ be a locally compact separable metrizable space.

\begin{definition} \label{defstrat} A $C^{\infty}$-stratification of $X$ is a pair
  $(\fS,N)$ such that:
\begin{enumerate}
 \item $\fS=\{s_i\}$ is a locally finite partition of
  $X$ into locally closed subsets of $X$, called the strata, which are
  smooth manifolds and which satisfies: 
  $$
    s_0\cap \bar{s_1}\not=\emptyset
    \mbox{ if and only if } s_0\subset \bar{s_1} .
  $$
  In that case we will write $s_{0}\leq s_{1}$ and $s_{0}<s_{1}$ if
  moreover $s_{0}\not= s_{1}$.  
\item $N=\{ \cN_s,\pi_s,\rho_s\}_{s\in \fS}$ is the set of control
  data or tube system: \\ 
  $\cN_s$ is an open neighborhood of $s$ in
  $X$, $\pi_s:\cN_s \rightarrow s$ is a continuous retraction and
  $\rho_s:\cN_s \rightarrow [0,+\infty[$ is a continuous map such that
  $s=\rho_s^{-1}(0)$. The map $\rho_s$ is either surjective or
  constant equal to $0$.\\ 
  Moreover if $\cN_{s_0}\cap s_1\not= \emptyset $ then the map 
  $$
    (\pi_{s_0},\rho_{s_0}):\cN_{s_0}\cap s_1
    \rightarrow s_0 \times ]0,+\infty[
  $$
  is a smooth proper submersion.
\item For any strata $s, t$ such that $s<t$, the inclusion
  $\pi_{t}(\cN_{s}\cap\cN_{t})\subset \cN_{s}$ is true and the
  equalities:
  $$
    \pi_{s}\circ \pi_{t}=\pi_{s} \makebox{ and } 
    \rho_{s}\circ\pi_{t}=\rho_{s}
  $$ 
  hold on $\cN_{s}\cap \cN_{t}$.
\item For any two strata $s_0$ and $s_1$ the following equivalences hold:
  $$
    s_0\cap \bar{s_1}\not= \emptyset \mbox{ if and only if }
    \cN_{s_0}\cap s_1 \not= \emptyset \ , 
  $$
  $$
    \cN_{s_0}\cap \cN_{s_1} \not= \emptyset \mbox{ if and only if }
   s_0\subset \bar{s_1},\ s_0=s_1 \mbox{ or } s_1\subset \bar{s_0}.
  $$
\end{enumerate}
\end{definition}

\medskip \noindent A stratification gives rise to a filtration: let $X_j$ be the union of strata of dimension
$\leq j$,  then:
 $$
  \emptyset \subset X_0 \subset \cdots \subset X_n =X \ .
 $$ 
We call $n$ the {\it dimension} of $X$ and $\smX:=X\setminus
X_{n-1}$ the {\it regular } part of $X$. The strata included in
$\smX$ are called {\it regular} while strata included in
$X\setminus \smx$ are called {\it singular}. The set of singular (resp. regular)
strata is denoted
$\fS_{sing}$ (resp.  $\fS_{reg}$).\\
For any subset $A$ of $X$, $\smA$ will denote $A\cap \smX$.

\medskip \noindent A crucial notion for our purpose will be the notion
of {\it depth}. Observe that the binary relation $s_{0}\leq s_{1}$ is
a partial ordering on $\fS$. 

\begin{definition}
The depth $d(s)$ of a stratum $s$ is the biggest
$k$ such that one can find $k$ different strata $s_0,\cdots,s_{k-1}$ such that
$$s_0<s_1<\cdots<s_{k-1}<s_k:=s .$$ 
The depth of the stratification $(\fS,N)$ of $X$ is:
$$d(X):= sup \{d(s),\ s\in \fS \} .$$
A stratum whose depth is $0$ will be called minimal.
\end{definition}

\smallskip \noindent We have followed the terminology of \cite{BHS},
but remark that the opposite convention for the depth also exists \cite{Ver}.

\medskip \noindent Finally we can define stratified pseudomanifolds:

\begin{definition}\label{defpseudomanifold}
A stratified pseudomanifold is a triple $(X,\fS,N)$ where $X$ is a
locally compact separable metrizable space, $(\fS,N)$ is a $C^{\infty}$-stratification
on $X$ and the regular part $\smx$ is a dense open subset of $X$.
\end{definition}

\noindent If $(X,\fS_X,N_X)$ and $(Y,\fS_Y,N_Y)$ are two stratified
pseudomanifolds  an homeomorphism $f:X\rightarrow Y$ is an isomorphism
of stratified pseudomanifold if: 
\begin{enumerate}\item
  $\fS_Y=\{f(s),\ s\in \fS_X\}$ and the
  restriction of $f$ to each stratum is a diffeomorphism onto its image.
\item $\pi_{f(s)}\circ f =f\circ \pi_s$ and $\rho_s=\rho_{f(s)}\circ
  f$ for any stratum $s$ of $X$.
\end{enumerate} 

\smallskip \noindent Let us make some basic remark on the previous definitions.

\begin{remark}\label{remdefstrat}\begin{enumerate}\item At a first sight, the definition of
    a stratification given here seems more restrictive than the usual
    one. 
In fact according to \cite{Ver} these definitions are equivalent.
\item Usually, for example in \cite{GoMa}, the extra assumption
  $X_{n-1}=X_{n-2}$ is required in the definition of stratified
  pseudomanifold. Our constructions remain without this extra assumption. 
\item A stratum $s$ is regular if and only if $\cN_s=s$ and then
  $\rho_s=0$.
\item  Pseudomanifolds of depth $0$ are smooth manifolds, and the
  strata are then union of connected components.
\end{enumerate}
\end{remark}

\noindent The following simple consequence of the axioms will be usefull enough
in the sequel to be pointed out: 
\begin{proposition}\label{orderedtubes}
  Let $(X,\fS,N)$ be  a stratified pseudomanifold. Any subset $\{s_i\}_I$ of distinct
  elements of $\fS$ is
  totally ordered by $<$ as soon as the intersection $\cap_{i\in
    I}\cN_{s_i}$ is non empty. In particular if the strata $s_0$ and $s_1$
  are such that $\cN_{s_0}\cap \cN_{s_1}\not= \emptyset$ then
  $d(s_0)\not=d(s_1)$ or $s_0=s_1$.
\end{proposition}

\medskip \noindent By a slight abuse of language we will
sometime
talk about a stratified 
pseudomanifold $X$ while we only have a partition $\fS$ on the space
$X$. This means that one can find at least one control data $N$ such
that $(X,\fS,N)$ is a stratified pseudomanifold in the sense of
our definition \ref{defpseudomanifold}.

\subsection{Examples}\label{exestrat} 

(1) Smooth manifolds are, without
  other mention, pseudomanifolds of depth $0$ and with a single stratum.\\
(2) Stratified pseudomanifolds of depth one are {\it wedges} and
    are obtained as follows.\\ Take $M$ to be a manifold with a
    compact boundary $L$ and let $\pi$ be a surjective submersion of
    $L$ onto a  manifold $s$. Consider the {\it mapping cone} of $(L,\pi)$~:
  $$c_{\pi}L:= L\times[0,1]/\sim_{\pi}$$ where $(z,t)\sim_\pi (z',t')$
  if and only if 
$(z,t)=(z',t')$ or 
$t=t'=0$ and $\pi(z)=\pi(z')$. The image of $L\times \{0\}$
  identifies with $s$ and by a slight abuse of notation we will denote
  it $s$. Now glue $c_{\pi}L$ and $M$ along their boundary in order to
  get $X$. The space $X$ with the partition $\{s,X\setminus s\}$ is a
  stratified pseudomanifold. \\ Two extreme examples are obtained by
  considering $\pi$ either equal to identity, with $s=L$ or equal to
  the projection on one point $c$. In the first case $X$ is a manifold
  with boundary $L$ isomorphic to $M$ and the stratification
  corresponds to
  the partition of $X$ by $\{L,X\setminus L\}$. In the second case $X$
  is a {\it conical manifold} and the stratification corresponds to 
  the partition of $X$ by $\{c,X\setminus c\}$, where $c$ is the
  singular point. \\
(3) Manifolds with
  corners with their partition into faces are stratified
  pseudomanifolds \cite{M1990,Mo1}. \\
(4) If $(X,\fS,N)$ is a pseudomanifold and $M$ is a smooth
  manifold then $X\times M$ is naturally endowed with a structure of
  pseudomanifold of same depth as $X$ whose strata are $\{s\times M,\
  s\in \fS\}$.\\
(5) If $(X,\fS,N)$ is a pseudomanifold of depth $k$
then $CrX:=X\times S^1 / X\times \{p\}$ is naturally endowed with a structure of
  pseudomanifold of depth $k+1$, whose strata are $\{s\times
  ]0,1[,\ s\in \fS\}\cup\{ [p] \}$. Here we have identified $S^1\setminus
  \{p\}$ with $]0,1[$ and we have denoted by $[p]$ the image of $
  X\times \{p\}$ in $CrX$.

\smallskip \noindent For example, if $X$ is the square we get the
following picture:

\medskip
\centerline{\includegraphics[width=11cm]{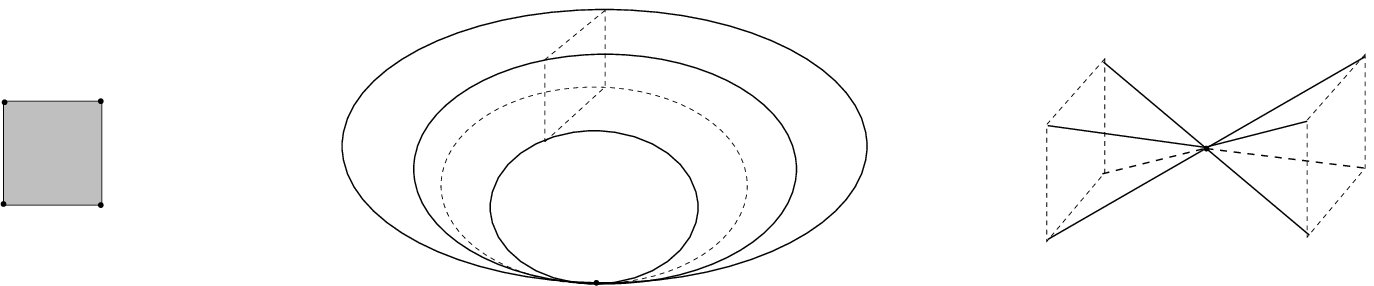}}

\subsection{The unfolding process}\label{unfolding}

Let $(X,\fS,N)$ be a stratified pseudomanifold. If $s$ is a singular
stratum, we let $L_s:=\rho_s^{-1}(1)$. Then $L_s$ inherits from $X$ a
structure of stratified pseudomanifold. \\ One can then define the
{\it open mapping
  cone} of $(L_s,\pi_s)$: 
$$ c_{\pi_s}L_s:=L_s\times [0,+\infty[ /
\sim_{\pi_s} $$ 
where $\sim_{\pi_s} $ is as before.

\smallskip \noindent According to \cite{Ver}, see also \cite{BHS}
the open mapping cone is naturally endowed  with a structure of stratified
pseudomanifold whose strata are $\{(t\cap L_s) \times ]0,+\infty[,\
t\in \fS\} \cup \{s\}$. 
Here we identify $s$ with the image of $L_s\times \{0\}$ in
$c_{\pi_s}L_s$.
Moreover,
up to isomorphism, the control data
on $X$ can be chosen such that one can find
a continuous retraction $f_s:\cN_s \setminus s \rightarrow L_s$ for which the map 
\begin{equation}\label{mapingcylindre} \begin{array}{cccc} \Psi_s:
    & \cN_s &  \rightarrow & c_{\pi_s}L_s   \\
  & z & \mapsto
  & \left\{ \begin{array}{cc} [f_s(z),\rho_s(z)] & \mbox{ if } z\notin s \\  z & \mbox{
    elsewhere }\end{array} \right. \end{array}\end{equation} 
is an isomorphism of stratified
pseudomanifolds. Here $[y,t]$ denotes the class in $c_{\pi_s}L_s$ of $(y,t)\in L_s\times
[0,+\infty[$.

\smallskip \noindent This result of local triviality  
around strata will be crucial for our
purpose. In particular it enables one to make the unfolding process
\cite{BHS} which consists in replacing each minimal stratum $s$ by $L_s$.
Precisely suppose that $d(X)=k>0$ and let $\fS_0$ be the set of
strata of depth $0$. Define $O_0:=\cup_{s\in \fS_0} \{z\in \cN_s \ \vert \ \rho_s(z) <
1\}$, $X_b=X\setminus O_0$ and 
$L:=\cup_{s\in \fS_0} \{z\in \cN_s \ \vert \ \rho_s(z) = 1 \}\subset X_b$. Notice that
it follows from remark \ref{remdefstrat} that the $L_s$'s where $s\in
\fS_0$ are disjoint and thus $L=\sqcup_{s\in
  \fS_0}L_{s}$. We let $$ 2X=X_b^- \cup L\times
[-1,1] \cup X_b^+ $$ 
where $X_b^\pm=X_b$ and $X_b^-$ (respectively $X_b^+$) is glued along $L$ with
$L\times \{-1\} \subset L\times [-1,1]$ (respectively $L\times \{1\}
\subset L\times [-1,1] $). \\
Let $s$ be a stratum of $X$ which is not minimal and which intersects $O_0$. We define the
following subset of $2X$: $$\tilde{s}:=
(s\cap X_b^-) \cup (s\cap L) \times [-1,1] \cup (s\cap X_b^+) $$
We then define $$\fS_{2X}:= \{ \tilde{s};\ s\in \fS \mbox{ and } s\cap
O_0\not= \emptyset\}\cup \{ s^-, s^+; s^\pm=s \in \fS \mbox{ and } s\cap
O_0= \emptyset\}.$$
The space $2X$ inherits from $X$ a structure of stratified
pseudomanifold of depth $k-1$ whose set of strata is $\fS_{2X}$.

\smallskip \noindent Notice that there is a natural map $p$ from $2X$
onto $X$. The restriction $p$ to any copy of $X_b$ is identity and
for $(z,t)\in L_s\times [-1,1]$, $p(z,t)=\Psi_s^{-1}([z,\mid \! t\! \mid
])$. The strata of $2X$ are the connected components of the pre-images by
$p$ of the strata of $X$. 

\medskip \noindent The interested reader can find all the details
related to the unfolding process in \cite{BHS} and \cite{Ver}
  where it is called decomposition. In particular starting
with a compact pseudomanifold $X$ of depth $k$, one can iterate this
process $k$ times and obtain a compact smooth
manifold $2^kX$ together with a continuous surjective map $\pi: 2^kX
\rightarrow X$ whose restriction to $\pi^{-1}(\smx)$ is a trivial
$2^k$-fold covering.

\begin{example}
Look at the square $C$ with stratification given by its  vertices, edges
and its interior.  It can be endowed with a
structure of stratified pseudomanifold of depth $2$. Applying
once the
unfolding process gives a sphere with $4$ holes: $S:=S^2\setminus
\{D_1,D_2,D_3,D_4\}$ where the $D_i$'s are disjoint and homeomorphic to open
disks. The set of strata of $S^2$ is then $\{\overset{\circ}{S},S_1,S_2,S_3,S_4\}$
where $S_i$ is the boundary of $D_i$ and $\overset{\circ}{S}$ the
interior of $S$. Applying the unfolding process
once more gives the torus with three holes.

\medskip
\centerline{\includegraphics[width=10cm]{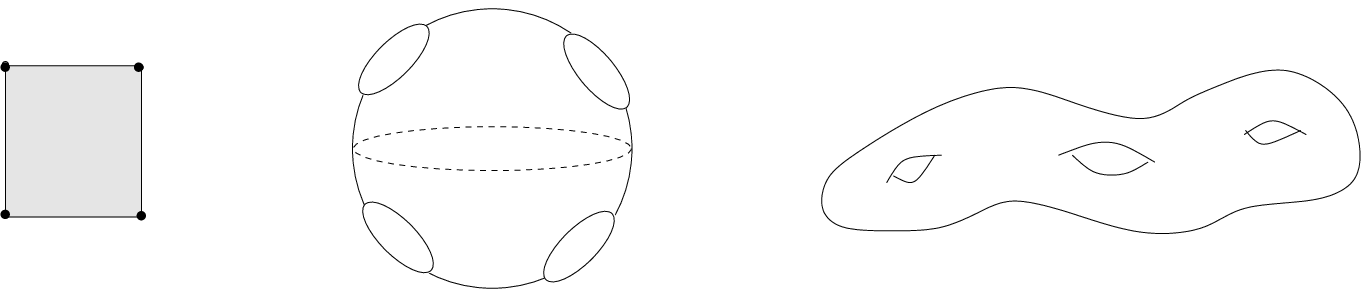}}
\end{example}



\section{The tangent groupoid and $\fS$-tangent space of a
 compact stratified pseudo-manifold}\label{tangentspace}

\subsection{The set construction}  We begin by the description at the
set level of the {\it $\fS$-tangent groupoid} and the {\it $\fS$-tangent space}
of a compact stratified pseudo-manifold.

\smallskip \noindent We keep the notation of the previous section: $X$ is a
compact stratified
pseudo-manifold, $\fS$ the set of strata, $\smx$ the regular part and
$N=\{\cN_s,\pi_s,\rho_s\}_{s\in \fS}$ the set of control data.

\smallskip \noindent
For each $s\in \fS$ we let $$O_{s}:=\{z\in \cN_{s} \ \vert \ \rho_{s}(z)<1\}
\mbox{ and } F_{s}:=O_{s} \setminus \bigcup_{s_0<s} O_{s_0} \ .$$

\noindent Note that $F_{s}=O_{s}$ if and only if $s$ is a minimal
stratum and $O_s=s$ when $s$ is regular. 

\begin{lemma} The set $\{F_s\}_{s\in \cS}$ form a partition of $X$.
\end{lemma}
\begin{proof} 
If $z$ belongs to $X$, let $\fR_z:=\{ s\in \fS\ \vert \ z\in \cN_s \
\mbox{ and }\ \rho_s(z)<1\}$.
It follows from proposition \ref{orderedtubes} that $\fR_z$ is a finite
set totally ordered by $<$. Since the set $\fR_z$ contains the stratum passing
through $z$, it is nonempty. Let $s_0^z$ be the minimal element of $\fR_z$.
Then $z$ belongs to $F_{s_0^z}$. Moreover, for all stratum $s\in \fS$, if
$s\not=s_0^z$ and $z\in O_s$, then $s\in \fR_z$, whence $s_0^z<s$, so that
$s\notin F_s$. 
\end{proof}

\noindent Recall that $ O_s^{\circ}=O_s \cap \smx$. We denote again by
$\pi_s: O_s^{\circ}\rightarrow s$ 
the projection. When $s$ is a
stratum, $\pi_s$ is a proper submersion and one can consider the
pull-back groupoid 
$\pb{\pi_s}(Ts) \rightrightarrows O_s^{\circ}$ of the usual tangent space
$Ts \rightrightarrows s$ by $\pi_s$. It is naturally endowed with a
 structure of smooth groupoid. When $s$ is a regular stratum, 
$s=O_s=O_s^{\circ}$ and $\pi_s$ is the identity map,  thus
$\pb{\pi_s}(Ts)\simeq TO_s^{\circ}$ in a canonical way. 

\smallskip \noindent At the set level, the $\fS$-{\it tangent space}
of $X$ is the groupoid: 
$$T^{\fS}X=\bigcup_{s\in \fS}
\pb{\pi_s}(TS)\vert_{F_s^{\circ}}\rightrightarrows \smx $$ 
where $F_s^{\circ}=F_s\cap \smx$. Following the cases of smooth
manifolds \cite{Co0} and isolated conical singularities
\cite{DL2}, the $\fS$-{\it tangent groupoid} of $X$ is defined to be a
deformation of the pair groupoid of the regular part of $X$ onto its $\fS$-tangent space:  
 $$
   \cG^t_X:= T^{\fS}X\times\{0\} \cup \smx \times \smx \times ]0,1]
   \rightrightarrows \smx \times [0,1].
 $$ 
 
\begin{examples} 
 \begin{enumerate}[\rm (1)] 
   \item When $X$ has depth
    $0$, we recover the usual tangent space and tangent groupoid.
   \item Suppose that $X$ is a trivial wedge (see example \ref{exestrat}):
     $$X=c_{\pi}L\cup M$$ 
    where $M$ is a manifold with boundary $L$ and $L$ is the
    product of two manifolds $L=s\times Q$ with
    $\pi:L\rightarrow s$ being the first projection. We have denoted by $c_{\pi}L=L\times
    [0,1]/\sim_{\pi}$ the mapping cone of $(L,\pi)$. In other word
    $c_{\pi}L=s\times cQ$ where $cQ:=Q\times[0,1]/Q\times\{0\}$ is the
    cone over $Q$. We denote again
    by $s$ the image of $L\times \{0\}$ in $X$. Then $X$ admits two
    strata: $s$ and $\smx=X\setminus s$, $F_s=O_s=L\times ]0,1[$ and
    $F_{\smX}=\smx\setminus O_s = M$. The tangent space is 
    $$
      T^{\fS}X=Ts\times (Q\times]0,1[) \times (Q\times]0,1[) \sqcup T M
      \rightrightarrows \smX
    $$ where $Ts\times (Q\times]0,1[) \times
    (Q\times]0,1[)$ is the product of the tangent space $Ts\rightrightarrows s$ with the pair
    groupoid over $Q\times ]0,1[$ and $TM$ denotes the restriction of the
    usual tangent bundle $T\smx$ to the sub-manifold with boundary $M$. 
\end{enumerate}
\end{examples}

\begin{remark} For any stratum $s$, the restriction of $\cG^t_X$ to 
  $\smF_s$ is equal to $$\pb{\pi_s}(TS)\vert_{F_s^{\circ}}\times \{0\}
   \cup \smF_s \times \smF_s \times ]0,1] \rightrightarrows \smF_s
  \times 
  [0,1]$$  
which is also the restriction to $\smF_s$ of $\pb{(\pi_s\times
  \id)}(\cG^t_s)$, the pull-back by 
  $\pi_s\times \id: O_s^{\circ}\times [0,1] \rightarrow S\times [0,1]$
  of the (usual) tangent groupoid of $s$:$$\cG^t_s=Ts\times \{0\}
  \cup s\times s 
  \times ]0,1] \rightrightarrows s\times [0,1]\ .$$  

\noindent In the following, we will denote  
  by $\cA^t_{\pi_s\times \id}$ the Lie algebroid of $\pb{(\pi_s\times
    \id)}(\cG^t_s)$. 
\end{remark} 

\subsection{The Recursive construction.} Thanks to the unfolding process
described in \ref{unfolding}, one can also construct the $\fS$-tangent
spaces of stratified pseudomanifolds by an induction on
the depth.

\smallskip \noindent If $X$ is of depth $0$, it is a smooth manifold
and the $\fS$-tangent space is the usual tangent space $TX$ viewed as
a groupoid on $X$.

\smallskip \noindent Let $k$ be an integer and assume that the $\fS$-tangent space of any
pseudomanifold of depth smaller than $k$ is defined. Let $X$ be  a stratified pseudomanifold of
depth $k+1$ and let $2X$ be the stratified pseudomanifold of
depth $k$ obtained from $X$ by applying \ref{unfolding}. With the notations of \ref{unfolding} we define 
$$T^{\fS}X=T^{\fS}2X\vert_{2\smX \cap X _b ^+} \underset{s\in S_0}\cup
\pb{\pi_s}(Ts)\vert_{\smO_s} \rightrightarrows \smX \ $$
where $T^{\fS}2X$ is the $\fS$-tangent space of the stratified
pseudomanifold $2X$. Here we have identified $2\smX \cap X
_b ^+$ with the subset $\smX\setminus O_0=X_b\cap \smX $ of $\smX$.
It is a simple exercise to see that this construction leads to the same definition of
$\fS$-tangent space as the previous one.

\subsection{The smooth structure}\label{subsection:smooth-structure} In this subsection we prove that 
the $\fS$-tangent space of a stratified pseudomanifold,  as well as
its $\fS$-tangent groupoid,  can be endowed with a smooth structure which reflects
the local structure of the pseudomanifold itself.    

Let $(X, \fS, N)$ be a stratified pseudomanifold.  
The smooth structure of $T^{\fS}X$ will depend on the stratification
and a smooth, decreasing, positive function $\tau:\RR \rightarrow
\RR$ such that $\tau([0,+\infty[)=[0,1]$, $\tau^{-1}(0)=[1,+\infty[$ and $\tau'$
does not vanish on $]0,1[$. The function $\tau$ will be called a {\it
  gluing function}.  We will also use functions associated with $\tau$
and defined on $\cN_{s}$ for any singular stratum $s$ by:
for each singular stratum:
 $$\tau_s=\tau \circ \rho_s $$
Observe that $\tau_s=0$ outside $\smO_s$.

\medskip \noindent Before coming into the details of the smooth
structure of $T^{\fS}X$,  let us describe its
consequences for  the convergence of sequences:\\ 
 A sequence $(x_n,V_n,y_n)\in \pb{\pi_{s_n}}(Ts_n)\vert_{F_{s_n}}$ where $n$ belongs to $\NN$, goes
 to $(x,V,y)\in \pb{\pi_{s}}(Ts)\vert_{F_{s}}$ if and only if:
\begin{equation}\label{top-TsX-sequences}
   x_n \rightarrow x,\ y_n \rightarrow y,\ V_n +
   \frac{\pi_s(x_n)-\pi_s(y_n)}{\tau_{s_n}(x_n)} \rightarrow V 
\end{equation}
The first two convergences have an obvious meaning, and they imply
that  for $n$ big enough, $s_{n}\le s$. The third one needs some
explanations.  Let us note
$z=\pi_{s}(x)=\pi_{s}(y)$ and
$z_{n}=\pi_{s_{n}}(x_{n})=\pi_{s_{n}}(y_{n})$. Since $\pi_s(x_n)$ and
$\pi_s(y_n)$ become close to $z$,  we can interpret
$w_{n}=\pi_s(x_n)-\pi_s(y_n)$
 as a vector in $T_{\pi_{s}(y_{n})}s$ (use any local chart of
$s$ around $z$). Moreover, using $\pi_{s_n}\circ \pi_s=\pi_{s_n}$,  
we see that this vector $w_{n}$ is vertical for $\pi_{s_{n}}$, that
is, belongs to the kernel ${\sf K}_n$ of the differential of
$\pi_{s_{n}}$ (suitably restricted to $s_{n}$). Now, the meaning of
last convergence in (\ref{top-TsX-sequences}) is 
$T\pi_{s_{n}}(V-w_{n}/\tau_{s_n}(x_n))-V_{n}\to 0$ which has to be
interpreted for each subsequences of $(x_n,V_n,y_n)_{n}$ with
$s_{n}=s_{n_{0}}$ for all $n\ge n_{0}$ big enough.

\medskip\noindent The smooth structure of $T^{\fS}X$  will be obtained by an induction
on the depth of the stratification,  and a concrete atlas will be given.  For the sake of
completeness,  we also explicit a Lie algebroid whose integration
gives the tangent groupoid $\cG^{t}_{X}$.  
We begin by describing the local structure of $\smx$ around its
strata,  then we will prove inductively the existence of a smooth
structure on the ${\fS}$-tangent space. Next,  an atlas of the
resulting smooth structure is given by  
brut computations.  A similar construction is easy to guess for
the tangent groupoid $\cG^{t}_{X}$. In the last part the previous
smooth structure is recovered in a more abstract approach using an
integrable Lie algebroid.

\medskip\noindent These parts are quite technical and can be
left out as soon as you believe that the tangent space and the tangent
groupoid can be endowed with a smooth structure compatible
with the topology described above.

\medskip\noindent Before going into the details, we should point out that the constructions described
above depend on the set of control data together with the choice of $\tau$, the
gluing function. As far as we know, there is no way to get rid
of these extra data. Nevertheless, a consequence of the last chapter is
that up to $K$-theory the $\fS$-tangent space $T^{\fS}X$ only depends
on $X$.
 
\subsubsection{The local structure of $\smx$.}\label{local-structure}

We now describe local charts of $\smx$ adapted to the
stratification,  called {\sl distinguished charts}. 

Let $z\in\smx$ and consider the set
 $$\fS_z:=\{s\in \fS \ \vert \ z\in \cN_s \ \mbox{ and }\ \rho_s(z) \leq 1 \}$$
It is a non empty  finite set, totally ordered according to
proposition \ref{orderedtubes},  thus we can write  
 $$ 
  \fS_z=\{s_0,\cdots,s_{\kappa} \}, \quad
   s_0<s_1<\cdots<s_{\kappa}$$ 
where $s_{\kappa}\subset \smx$ must be regular. Let $n_{i}$ be the
dimension of $s_{i}$,  $i\in\{0, 1,\ldots,\kappa\}$ and $n=n_{\kappa}=\dim\smx$. 

Let $U_z$ be an
open neighborhood of $z$ in $\smx$ such that the following hold:
\begin{equation}
  \label{special-neigh}
   U_z\subset \bigcap_{s\in\fS_{z}} \cN_{s}\qquad \text{ and
   }\qquad 
   \forall s\in\fS_{\mathrm{sing}}, \  U_z\cap O_s \not=
   \emptyset\Leftrightarrow s\in \fS_z
\end{equation}
In particular,  the following hold on $U_{z}$:
\begin{equation}
  \label{compcond}  \mbox{ for } 0\leq i\leq j \leq {\kappa}:  \pi_{s_i}\circ \pi_{s_j}=\pi_{s_i} \makebox{
    and } \rho_{s_i}\circ 
  \pi_{s_j}=\rho_{s_i} \ .
\end{equation}
Without loss of generality,  we can also assume that $U_{z}$ is the
domain of a local chart of $\smx$. 

\medskip \noindent If $\kappa=0$,  any local chart of $\smx$ with
domain $U_{z}$ will be called distinguished. When ${\kappa}\geq 1$,
we can take successively canonical forms of
the submersions $\pi_{s_0}, \pi_{s_1}\ldots, \pi_{s_{\kappa}}$
available on a possibly smaller
$U_z$, that is,  one can shrink $U_{z}$ enough and find diffeomorphisms:
\begin{equation}
  \label{disting-maps}
 \phi_{i}:\pi_{s_i}(U_z)\to\RR^{n_{i}} \text{ for all } i\in\{0, 1, \ldots, \kappa \}
\end{equation}
such that the diagram:
\begin{equation}\label{redressement-Op}\hspace{-1cm}
  \xymatrix{
   \pi_{s_i}(U_z)\ar[d]_{\pi_{s_j}}  \ar[r]^{\phi_{i}} & \RR^{n_{i}}\ar[d]^{\sigma_{n_{j}}} \\
    \pi_{s_j}(U_z) \ar[r]^{\phi_{j}}  & \RR^{n_{j}}   }
\end{equation}  
commutes for all $i, j \in\{0, 1, \ldots, \kappa \}$ such that $i\ge
j$.  Above, for any integers $p\ge d$,  the map $\sigma_d:\RR^p\to\RR^d$ denotes the canonical
projection onto the last $d$ coordinates.

Remember that $s_{\kappa}$ is regular so $\pi_{s_{\kappa}}$ is the
identity map and $\phi:=\phi_{\kappa}$ is a local chart around $z$ of
$\smx$. Now we set:
\begin{definition}
  A {\it distinguished chart} of $\smx$ around $z\in \smx$ is a local chart
  $(U_{z}, \phi)$ around $z$ such that $U_{z}$ satisfies
  (\ref{special-neigh}) together with diffeomorphisms (\ref{disting-maps})
  satisfying (\ref{redressement-Op}) and $\phi=\phi_{\kappa}$. 
\end{definition}
\smallskip \noindent From now on, a riemannian metric is chosen on
$\smx$ (any adapted metric in the sense of \cite{BHS1992} is suitable
for our purpose). 
Recall that for any stratum $s$, the map $\pi_s:\smN_s \rightarrow s$ is a 
smooth submersion. Thus,  if $K_s\subset T\smN_s$ denotes the kernel of the
differential map $T\pi_s$ and $q_s:T\smN_s \rightarrow T\smN_s$ the orthogonal
projection on $K_s$,  the map 
\begin{equation}\label{iso-tsks} (q_s, T\pi_s): T\smN_s \rightarrow
 K_s \oplus \pi_s^*(Ts)   \end{equation} is an isomorphism and 
the vector bundle $\pi_s^*(Ts)$ can be identified with the othogonal
complement of $K_{s}$ into $T\cN_{s}=T\smx|_{\cN_{s}}$.

Now,  let $(U_{z}, \phi)$ be a distinguished chart
around some $z\in\smx$. Set $\fS_{z}=\{s_{0}, s_{1}, \ldots,
s_{\kappa}\}$ with $s_{0}<s_{1}<\ldots<s_{\kappa}$,  and set
$K_{i}=K_{s_{i}}|_{U_{z}}$,  $U_{i}=\pi_{s_{i}}(U_{z})$ for all $i=0, 1, \ldots, \kappa$.   

By \ref{compcond} we have:
\begin{equation}  \label{tangent-filtration}
U_{z}\times\{0\}=K_{\kappa}\subset
  K_{\kappa-1}\subset \cdots \subset K_{1} \subset K_{0}\subset 
  TU_z. 
\end{equation} 
Rewriting the diagram (\ref{redressement-Op}) for the differential
maps and $i=\kappa$,  we get for all $j\le\kappa$:
\begin{equation}\label{redressement-Op-diff}\hspace{-1cm}
  \xymatrix{
   TU_z\ar[d]_{T\pi_{s_j}}  \ar[r]^-{T\phi} & \RR^{n}\ar[d]^{\sigma_{n_{j}}\times\sigma_{n_{j}}} \times\RR^{n} \\
    TU_j \ar[r]^-{T\phi_{j}}  & \RR^{n_{j}}\times \RR^{n_{j}}  }
\end{equation}  
and we see that $T\phi$ sends the filtration
(\ref{tangent-filtration}) to the following filtration:
\begin{equation}
  \label{euclidean-filtration}
   \RR^n\times\{0\}\subset \RR^n\times\RR^{n-n_{\kappa-1}}\subset\cdots\subset \RR^n\times\RR^{n-n_{0}}\subset \RR^n\times\RR^n, 
\end{equation}
where $\RR^{n-n_{i}}$ is included in $\RR^{n}$ by the map $v\mapsto (v, 0)\in\RR^{n-n_{i}}\times\RR^{n_{i}}\simeq\RR^{n}$. 
This property can be reformulated in terms of natural graduations
associated with (\ref{tangent-filtration}) and
(\ref{euclidean-filtration}) (and will be used in this latter
form). Indeed,  let   $T^{i}$  be the orthogonal complement of $K_{i}$ into $K_{i-1}$ for all
$i=0,\ldots,\kappa$ (with the convention $K_{-1}=TU_z$). Moreover,  on
the euclidean side, let us embed
$\RR^{n_{i}-n_{i-1}}$ into $\RR^{n}$ by the map:
 $$
  v\in \RR^{n_{i}-n_{i-1}}\longmapsto (0, v, 0)\in \RR^{n-n_{i}}\times\RR^{n_{i}-n_{i-1}}\times\RR^{n_{i-1}}\simeq\RR^{n}
 $$
for all $i=0,1, \ldots,\kappa$ (by convention $n_{-1}=0$).  With these
notations and conventions,  the filtrations
(\ref{tangent-filtration}) and (\ref{euclidean-filtration}) give rise
to the following decompositions: 
\begin{equation}
  \label{tangent-graduation}
  TU_z = T^{\kappa}\oplus T^{\kappa-1}\oplus \cdots\oplus
  T^{0} 
\end{equation} 
and 
\begin{equation}\label{Euclideangraduation}
  \RR^{n}\times\RR^n =\RR^{n}\times(\RR^{n-n_\kappa}\oplus \RR^{n_\kappa-n_{\kappa-1}}\oplus \cdots\oplus
   \RR^{n_1-n_0}\oplus \RR^{n_0})
\end{equation}
Now,  that $T\phi$ respects the filtrations
(\ref{tangent-filtration}) and (\ref{euclidean-filtration}) means
that for all $x\in U_{z}$ the linear map $T\phi_{x}$ is upper
triangular with respect to the decompositions
(\ref{tangent-graduation}) and (\ref{Euclideangraduation}). 

The {\sl diagonal blocks} of $T\phi$ are the maps:
\begin{equation}
  \label{diagonal-blocks}
  \delta^j\phi: T^{j} \longrightarrow \RR^n\times\RR^{n_{j}-n_{j-1}};
  \ j=0, 1,\ldots,\kappa, 
\end{equation}
obtained by composing $T\phi$ on the left and on the right respectively by
the projections:
 $$
   TU_z = T^{\kappa}\oplus T^{\kappa-1}\oplus \cdots\oplus
  T^{0}\longrightarrow T^{j}
 $$
and
 $$
 \RR^{n}\times(\RR^{n-n_\kappa}\oplus \RR^{n_\kappa-n_{\kappa-1}}\oplus \cdots\oplus
   \RR^{n_1-n_0}\oplus \RR^{n_0}) \longrightarrow \RR^{n}\times \RR^{n_{j}-n_{j-1}}. 
 $$
The {\sl diagonal part} of $T\phi$ will be defined by 
$\Delta\phi=(\delta^\kappa\phi,\delta^{\kappa-1}\phi,\ldots,\delta^0\phi)$.
Of course,  the inverse of $T\phi$ is also upper triangular with
diagonal blocks given by 
 $(\delta^j\phi)^{-1}$, $j=0,1,\ldots,{\kappa}$.

\medskip \noindent We have similar properties for all the underlying maps $\phi_i$, $i=0,1,\ldots,\kappa-1$ coming with 
the distinguished chart.  To fix notations and for future references,
let $U_i$ denote $\pi_{s_i}(U_z)$, and $T_{i}^{j}$ denote
$T\pi_{s_i}(T^{j})$ for all $j\le i<\kappa$. Applying now $T\pi_{s_i}$ to (\ref{tangent-graduation}) yields:
\begin{equation}
  \label{tangent-graduation-si}
  TU_i = T^{i}_{i}\oplus T_{i}^{i-1}\oplus \cdots\oplus T_{i}^{0},
\end{equation}
It follows that the differential maps: 
\begin{equation}\label{diff-phi-si}
 T\phi_{i}:T^{i}_{i}\oplus T_{i}^{i-1}\oplus
 \cdots\oplus T_{i}^{0} 
  \longrightarrow \RR^{n_{i}}\times (\RR^{n_{i}-n_{i-1}}\oplus \cdots\oplus
   \RR^{n_{1}-n_{0}}\oplus \RR^{n_{0}})
\end{equation}
for  all $i=0, 1, \ldots, \kappa-1$ are upper triangular with
diagonal blocks $\delta^j\phi_i$ defined as above. 
Note that for all  $j\le i\le k\le \kappa$,
$(T\pi_{s_{i}})(T^{j}_{k})=T^{j}_{i}$ and that applying the correct
restrictions and projections in (\ref{redressement-Op-diff}) gives the following
commutative diagram:
\begin{equation}\label{redressement-Op-delta}\hspace{-1cm}
  \xymatrix{
    T^{j}\ar[d]_{T\pi_{s_i}} \ar[r]^-{{\delta^{j}\phi_{k}}\; } &
    \RR^{n_{k}}\times \ar[d]^{\sigma_{n_{i}}\times\id} \RR^{n_{j}-n_{j-1}}  \\
    T^{j}_{i} \ar[r]^-{{\delta^{j}\phi_{i}}\; } &
   \RR^{n_{i}}\times \RR^{n_{j}-n_{j-1}}  }
\end{equation}  

\subsubsection{The smooth structure by induction}

We show that $T^{\fS}X$ can be provided with a smooth
structure by a simple recursive argument. 

\smallskip \noindent Let us first introduce the {\it $s$-exponential
  maps}. Let $s$ be a stratum. The corresponding {\it $s$-exponential
  map} will be an exponential along the fibers of $\pi_s$. Precisely,
recall that the map $\pi_s:\smN_s \rightarrow s$ is a  
smooth submersion,  $K_s\subset T\smN_s$ denotes the kernel of the
differential map $T\pi_s$ and $q_s:T\smN_s \rightarrow T\smN_s$ the orthogonal
projection on $K_s$. The subbundle $K_s$ of $T\smX$ inherits from
$T\smX$ a riemannian metric whose associated riemannian connection is
 $\nabla^s=q_s\circ \nabla$, where $\nabla$ is the riemannian
connection of the metric on $\smX$.
The associated exponential map  $$Exp^s:V_s\subset K_s \rightarrow
\smN_s\ $$ is smooth and defined on an open neighborhood $V_s$ of the zero section of
$K_s$. Moreover it satisfies:
\begin{itemize}\item[-] $\pi_s\circ
  Exp^s=\pi_s$.
\item[-] For any fiber $L^s$ of $\pi_s$, the
  restriction of $Exp^s$ to $L^s$  is the
  usual exponential map for the submanifold $L^s$ of $\smX$ with the
  induced riemannian structure.
\end{itemize}

\smallskip \noindent If $X$ is a stratified pseudomanifold of depth $0$ it is smooth
  and its $\fS$-tangent space is the usual tangent space 
$TX$ equipped with its usual smooth structure.

\smallskip \noindent Suppose that the $\fS$-tangent space of any stratified pseudomanifold
of depth strictly smaller than $k$ is equipped
with a smooth structure for some integer $k>0$ . Let $X$ be  a stratified pseudomanifold of
depth $k$ and take $2X$ be the stratified pseudomanifold of
depth $k-1$ obtained from $X$ by the unfolding process
\ref{unfolding}. According to \ref{tangentspace}, with the notations of \ref{unfolding} we have
$$T^{\fS}X=T^{\fS}2X\vert_{2\smX \cap X _b ^+} \underset{s\in \fS_0}\cup
\pb{\pi_s}(Ts)\vert_{\smO_s} \rightrightarrows \smX \ .$$
Let $\sml$ be the boundary of $2\smX \cap X_b^+$ in $\smX$. We
equip the restriction of $T^{\fS}X$ to $2\smX \cap X _b ^+
\setminus \sml$ with the smooth structure coming from $T^{\fS} 2X$ and its
restriction to any $O_{s_0}$, $s_0\in \fS_0$,
with the usual smooth structure. It remains to describe the gluing over $\sml$. 
One can find an open subset $W$ of $T^{\fS} 2X$ which contains the
restriction of $T^{\fS} 2X$ to  $\sml$  
such that the following map is defined:
$$
  \begin{array}{cccc} \Theta: & W &
  \longrightarrow & T^{\fS}X \\ & (x,u,y) &
  \mapsto & \left\{ \begin{array}{ll}
      (x,T\pi_{s_0}(u),Exp^{s_0}(y,-\tau_{s_0}(x)q_{s_0}(u))) \mbox{ if
        }
        x \in \smO_{s_0},\  s_0\in \fS_0 \\ (x,u,y) \mbox{
          elsewhere }\end{array}
    \right. 
  \end{array}
$$ 
Here, if $s$ denotes the unique stratum such that $x, y\in F_{s}$,
the vector bundle $\pi_s^*(Ts)$ is identified with the orthogonal complement of $K_s$ into
$T\smN_s$, in other words $q_{s_0}(y,u)=q_{s_0}(W-q_s(W))$ where
$W\in T_y\smX$ satisfies $T\pi_s(W)=u$. 

Then, we equip $T^{\fS}X$ with the 
unique smooth structure compatible with the one previously defined
on $T^{\fS}X|_{\smX\setminus\sml}$
and  such that the map $\Theta$ is a smooth diffeomorphism onto its image.
The non trivial point is to check that the restriction of the map
$\Theta$ over $\smO_{s_0}$ is a diffeomorphism onto its image for any
$s_0\in \fS_0$. This will follows from the following
lemma. 

\begin{lemma} If $s_0<s$, for any $x_0\in
s_0$ and $x\in s$ with $\pi_{s_0}(x)=x_0$. The following assertions
hold:
\begin{enumerate}\item $E:=q_{s_0}(\pi_s^*(Ts))$ is a sub-bundle of
  $K_{s_0}$ of dimension $\dim(s)-\dim(s_0)$.
\item Let $E^x:=q_{s_0}(\pi_s^*(Ts))\vert_{\pi_s^{-1}(x)} $ be the restriction of $E$ to the submanifold
  $\pi_s^{-1}(x)$. There exists a  neighborhood $W$ of the zero section of
  $E^x$ such that the restriction of $Exp^{s_0}$ to $W$ is a
  diffeomorphism onto a neighborhood of $\pi_s^{-1}(x)$ in
  $\pi_{s_0}^{-1}(x_0)$.
\end{enumerate}
\end{lemma}

\begin{proof} 1. The first assertion follows from the inclusion: $\pi_{s_0}^*(Ts_0)
  = K_{s_0}^{\bot} \subset K_s^{\bot}=\pi_{s}^*(Ts)$ which ensures
  that the dimension of the fibers of $q_{s_0}(\pi_s^*(Ts))$ is
  constant equal to $\dim(s)-\dim(s_0)$.\\
The same argument shows that $K_{s_0}=K_s \oplus E$.\\
2. If $\Psi$ denotes the restriction of $Exp^{s_0}$ to $E^x$ then
$T\Psi(z,0)(U,V)=U+V$ where $(z,U)\in K_s$ and $V\in E_z$. Since
$K_s\cap E$ is the trivial bundle we get that $T\Psi$ is injective and
since $E^{x}$ and $\pi_{s_{0}}^{-1}(x_{0})$ have same dimension,  it
is  bijective. We conlude with the local
inversion theorem.
\end{proof}

\subsubsection{An atlas for $T^\fS X$} 

The atlas will contain  two kinds of local charts. The kind of these
charts will depend on
the fact that their domains meet or not a gluing between the different
pieces composing the tangent space $T^\fS X$, that is the boundary of
some $F_s$. 

\smallskip \noindent The first kind of charts, called {\sl regular
  charts}  are charts whose domain is contained in  
$T^\fS X|_{\interior{F_s}}$ for a given stratum $s$ of the stratification. We
observe that $T^\fS X|_{\interior{F_s}}$ is a smooth groupoid as an
open subgroupoid of $\pb{\pi_s}(Ts)\rightrightarrows \smN_s$.
Thus, regular charts have domains  contained in 
 $$ 
    \sqcup_{s\in\cS}\interior{F_s}
 $$
and coincide with the usual local charts of the (disjoint) union of the smooth
groupoids $\sqcup_{s\in\cS}\pb{\pi_s}(Ts)$. 
 
\smallskip \noindent The second kind of charts, called {\sl
  deformation charts} (adapted to 
a stratum $s$), are charts
whose domain meets $T^\fS X|_{\partial F_s }$ for a given
stratum $s$, that is, charts around points in 
$$
\bigcup_{s\in\fS} T^\fS X|_{\partial F_s }. 
$$
Their description is more involved. Let $(p,u,q)\in T^\fS X$. Thus
there is a stratum $s$ such that $p$ and $q$ belong to $F_s$ with
$\pi_s(p)=\pi_s(q)$ and $u\in T_{\pi_s}(p)s$. Assume that $p\in
\partial F_{s}$.  This means that $\rho_{s}(p)<1$,  that $\rho_{t}(p)\ge 1$ for
all strata $t<s$ and that the set of strata $t$ such that $t<s$ and $\rho_{t}(p)=1$ is not empty. Using again the axioms of the stratification, we see that this set is totally ordered and we denote $s_0,s_1,\ldots,s_{l-1}$ its elements listed by increasing order. We also set $s_l=s$. Observe that:
 \begin{equation} \label{chart-filtration}
 \{ s_0,s_1,\ldots, s_l\} = \fS_p\cap \{t\in\fS\ | \ t\le s\} 
\end{equation}
and that, thanks to the compatibility conditions (\ref{compcond}),
this set only depends on $\pi_s(p)$ and thus is equal with the corresponding set associated with $q$. 

\smallskip \noindent Let us take  distinguished charts $\phi: U_p
\rightarrow \RR^n$ around $p$ and $\phi': U_q \rightarrow \RR^n$
around $q$. Since $\pi_s(p)=\pi_s(q)$, we can also assume without loss of generality that:
\begin{equation}\label{factorisation} 
   \pi_{s_i}(U_p)=\pi_{s_i}(U_q) \hbox{ and }
   \phi_{i}=\phi'_{i} \hbox{ for } i=0,\ldots,l.
\end{equation}

\noindent We will use the same notations as in paragraph \ref{local-structure}: $n_i=\dim s_i$, $U_i= \pi_{s_i}(U_p)$, 
$K_i=\ker(T\pi_{s_i})|_{U_p}$, $T^i=K_i^{\perp K_{i-1}}$ for all
$i=0,1,\ldots,l$ (here again $K_{-1}=TU_p$). The main difference with
the settings of the paragraph \ref{local-structure}is that we forget the strata bigger than $s$ in $\fS_p$ and  $\fS_q$  to concentrate on the lower (and common) strata in $\fS_p$ and  $\fS_q$. It amounts to forget the tail of the filtration (\ref{tangent-filtration}) up to the term $K_l$:
\begin{equation}
  \label{tangent-filtration2}
   K_l \subset K_{l-1}\subset \cdots\subset K_0\subset TU_p
\end{equation}
and this leads to a less fine graduation: 
\begin{equation}
  \label{tangent-graduation2}
  TU_p = K_l\oplus T^l\oplus T^{l-1}\oplus \cdots\oplus T^0
\end{equation}

\medskip \noindent Let us also introduce the positive smooth functions: 
 $$
   t_i  = \sum_{j=0}^i \tau\circ\rho_{s_j}\ , \  i=0,1,\ldots,l \ ;
   \qquad  \theta_i=\prod_{j=i-1}^{l}
   t_j \ , \ i=1,\ldots,l
 $$
Note that  $t_j$ (resp. $\theta_j$) is strictly positive on $F_{s_i}$
if $j\ge i$ (resp. $j>i$) and
vanishes identically if $j<i$ (resp. $j\le i$).

Finally we will write:
 $$
  \forall x\in U_{p}, \quad \phi(x)=(x^{l+1},x^l,\ldots,x^1,x^0)\in\RR^{n-n_{l}}\times\RR^{n_{l}-n_{l-1}}\times\cdots
   \times\RR^{n_{1}-n_{0}}\times\RR^{n_{0}}, 
 $$
and for all $j=0, 1, \ldots, l$, 
 $$
  \pi_{s_{j}}(x)=x_{j}, 
 $$
thus $\phi_{j}(x_{j})=(x^{j}, x^{j-1}, \ldots, x^{0})\in\RR^{n_{j}}$;
and we adopt similar notations for $\phi'$ and $y\in U_{q}$.

\medskip \noindent We are ready to define a deformation chart around the point $(p,u,q)$. The domain will be: 
\begin{equation}
  \label{domain-deformation-chart}
  \wtU = T^\fS X|_{U_q}^{U_p}
\end{equation}
and the chart itself:
\begin{equation}
  \label{notation-deformation-chart}
  \wtphi: \wtU \to \RR^{2n}
\end{equation}
is defined as follows. Up to a shrinking of $U_p$ and $U_q$, the following is true: for all $(x,v,y)\in \wtU$, there exists a
unique $i\in\{0,1,\ldots,l\}$ such that $x\in F_{s_i}$. Then $(x,v)\in \pi_{s_i}^*(TU_i)$, and we set:
\begin{equation}
  \label{definition-deformation-chart}
   \wtphi(x,v,y) =
   \left(\phi(x),\frac{x^{l+1}-y^{l+1}}{\theta_{l+1}(x)},\ldots,
   \frac{x^{i+1}-y^{i+1}}{\theta_{i+1}(x)},\Delta\phi_i(x_{i},v)\right) 
\end{equation}
The map $\wtphi$ is clearly injective with inverse defined as
follows. For  $(\bfx,\bfw)\in \wtphi(\wtU)$ and $i$ such that $\phi^{-1}(\bfx)\in F_{s_i}$:
\begin{eqnarray*}
  \wtphi^{-1}(\bfx,\bfw) =
   \left(\phi^{-1}(\bfx),(\Delta\phi_{i})^{-1}(\bfx_{i}, \bfw),
     \phi'^{-1}(\bfx-\Theta^{[i+1]}(\phi^{-1}(\bfx))\cdot \bfw)\right) 
\end{eqnarray*}
where $\bfx_{i}=\sigma_{n_{i}}(\bfx)$ and,  using the decomposition 
 $$
   \bfw=(\bfw^{l+1},\bfw^{l}\cdots,\bfw^{0})\in\RR^{n-n_{l}}\times\RR^{n_{l}-n_{l-1}}\times\cdots
   \times\RR^{n_{1}-n_{0}}\times\RR^{n_{0}}, 
 $$
we have set
$$
 \Theta^{[i+1]}(x)\cdot \bfw=
 \theta_{l+1}(x)\bfw^{l+1}+\cdots+\theta_{i+1}(x)\bfw^{i+1}\in
 \RR^{n-n_i}\times\{0\}\subset\RR^n.
$$
To ensure that $( \wtphi,\wtU)$ is a local chart, it remains to check
that $\wtphi(\wtU)$ is an open subset of $\RR^{2n}$.  
It is easy to see that $\wtphi(\interior{F_{s_i}})$ is open for every
$i\in\{0,\ldots,l\}$ so we consider  $(p,u,q)\in \wtU$ such that 
$p\in \pa F_{s_i} $ for some integer $i$. Let 
$J=\{i_0,\ldots,i_k\}\subset\{0,1,\ldots,i-1\}$ such that: 
 $$
   \forall j\in J,\quad \rho_{s_j}(p)=1\ .
 $$
Thus we have: 
\begin{equation}
  \label{ineq-eq-boundary-Fsi}
   \rho_{s_i}(p)<1; \quad  \forall j\in J,\ \rho_{s_j}(p)=1; \quad \forall j\not\in J 
  \hbox{ and  } j< i,  \ \rho_{s_j}(p)>1 
\end{equation}
by construction, $q$ satisfies the same relations.  Set
$\wtphi(p,u,q)=(\bfx_0,\bfv_0)$. Using the Taylor formula and the fact
that $\theta_{j+1}$ is negligible with respect to $1-\rho_{s_j}$ at the
region $\rho_{s_j}=1$, noting also the invariance of $\rho_{s_k}$ with respect
to perturbations of points along the fibers of
$\pi_{s_{k+1}},\pi_{s_{k+2}},\ldots$; we prove that there exist an open ball $B_1$ of
$\RR^n$ centered at $\bfx_0$ and an open ball $B_2$ of $\RR^n$
centered at $0$ and containing $\bfv_0$ such that for all 
$(\bfx,\bfv)\in B_1\times B_2$, if 
$$ 
  x = \phi^{-1}(\bfx) \in F_{s_j} \mbox{ for } j\in J \mbox{ or }
  j=i, \mbox{ then } 
 y= \phi'^{-1}(\bfx-\Theta^{[j+1]}(x)\cdot \bfv)\in F_{s_j}\ .
$$
This proves that $(\bfx,\bfv)\in\im \wtphi$, thus 
 $$
   \wtphi(p,u,q)\in B_1\times B_2\subset \im\wtphi
 $$
and the required assertion is proved. We end with:
\begin{theorem}
  The collection of regular and deformation charts provides $T^\fS X$
  with a structure of smooth groupoid. 
\end{theorem}
\begin{proof}
The compatibility between a regular and a deformation chart contains
no issue and is ommitted. We need only to check the compatibility 
between a deformation chart adapted to a stratum $s$ and a deformation
chart adapted to a stratum $t$, when their domains overlap, which
implies automatically that $s<t$ or $s>t$ or $s=t$. 

\noindent Let us work out only the case $s=t$, since the other case is similar. 
We have here to compare two  charts $\wtphi$
and $\wtpsi$ with common domain $\wtU$ and involving the same chain of
strata $s=s_l>s_{l-1}>\cdots>s_0$. The whole notations are as before
and $\psi,\psi'$ are the underlying charts of $X^\circ$ allowing the
definition of $\wtpsi$.  We note, for
the sake of concision, $u^k$ (resp. $u'^k$), $k=l+1,\ldots,0$, the coordinate functions of
$u:=\psi\circ\phi^{-1}$ (resp. $\psi'\circ\phi'^{-1}$) with respect to the decomposition
(\ref{Euclideangraduation}) of $\RR^n$. Observe, thanks to the
particular assumptions made on $\phi,\phi',\psi,\psi'$  
(cf.(\ref{redressement-Op}),
(\ref{factorisation})), that $u^k(\bfx)$ only depends on
$\bfx_k:=(\bfx^k,\bfx^{k-1},\ldots,\bfx^0)\in\RR^{n_{k}}$ and that $u^k=u'^k$ for all
$k<l+1$.  Let $(\bfx,\bfv)\in\im\wtphi$ and $i$ such
that $x=\phi^{-1}(\bfx)\in F_{s_i}$. Then: 
\begin{equation}\label{change-coordinates}
 \begin{matrix}
\wtpsi\circ\wtphi^{-1}(\bfx,\bfv)&  
  = \left(u(\bfx),
 \frac{u^{l+1}(\bfx)-u'^{l+1}(\bfx-\Theta^{[i+1]}\cdot\bfv)}{\theta^{l+1}},
  \frac{u^{l}(\bfx)-u^{l}(\bfx-\Theta^{[i+1]}\cdot\bfv)}{\theta^{l}},\ldots \right.\\
  & \left. \ldots,\frac{u^{i+1}(\bfx)-u^{i+1}(\bfx-\Theta^{[i+1]}\cdot\bfv)}{\theta^{i+1}},
 (\Delta\psi_{i})\circ (\Delta\phi_{i})^{-1}(\bfv)\right)
 \end{matrix}
\end{equation}
We need to check that  the above expression matches smoothly with the
corresponding expression for an integer $k\in[i,l]$ when $\theta_k(x)$ (and
thus $\theta_{k-1},\ldots,\theta_{i+1}$) goes to zero. For that, the
Taylor formula applied to $u^r$, $k\ge r\ge i+1$, shows that the map
defined below is smooth in $(\bfx,\bfv,t)$ where $(\bfx,\bfv)$ are as
before and $t=(t_l,t_{l-1},\ldots,t_0)\in\RR^{l+1}$ is this time an
arbitrary $(l+1)$-uple close to $0$:  
 $$  
   \begin{cases} 
     \frac{u^{r}(\bfx)-u^{r}(\bfx-\Theta^{[i+1]}\cdot\bfv)}{\theta^{r}} 
                          & \hbox{ if } \theta_r=\Pi_{r-1}^lt_j\not=0\\
     d(u^r)_{\bfx}(\bfv^r+t_{r-2}\bfv^{r-1}+\cdots+t_i\bfv^{i+1}) &
     \hbox{ if } \exists j\in\{r-1,r,\ldots,l\} \hbox{ such that } t_j=0.
   \end{cases}  
 $$
In our case, $t_j=t_j(x)$ and $t_{k-1},\ldots,t_{i}$ go to zero, so
the second line in the previous expression is just:
 $$ d(u^r)_{\bfx}(\bfv^r) $$
and for obvious matricial reasons:
 $$ 
   d(u^r)_{\bfx}(\bfv^r)
   =(\Delta\psi_{k})\circ(\Delta\phi_{k})^{-1}(\bfv^r)
 $$ 
Summing up these relations for $r=i+1,\ldots,k$, we arrive at the
desired identity. 

\smallskip \noindent Thus, $T^\fS X$ is endowed with a structure of smooth
manifold. Changing the riemannian metric on $\smx$ modifies the
choices of the $T^{i}_{j}$'s,  but gives rise to compatible charts. Moreover, the smoothness of all
algebraic operations associated with this groupoid is easy to check
in these local charts.   
\end{proof}

\subsubsection{The Lie algebroid of the tangent groupoid} 

We describe here the smooth structure of the tangent space
via its infinitesimal structure, namely its Lie algebroid. Precisely, we define 
$$\begin{array}{cccc} Q_s: & T\smx & \longrightarrow & T\smx  \\ 
 & (z,V) & \mapsto & \left\{ \begin{array}{ll} (z,\tau_s(z)q_s(z,V)) & \mbox{ if
   } z\in \smN_s \\  0 & \mbox{ elsewhere } \end{array}  \right. \end{array} $$
 
\noindent By a slight abuse of notation, we will keep the notations $q_s$ and $Q_s$
for the corresponding maps induced on the set of local tangent vector fields on
$\smx$.

\smallskip \noindent Let $\cA$ be the smooth vector bundle $\cA:=
T\smx \times [0,1]$ over 
$\smx \times [0,1]$. We define the following morphism of vector bundle
: 
$$\begin{array}{cccc} \Phi: & \cA=T\smx\times [0,1] &\longrightarrow 
  &T\smx \times T [0,1] \\ & (z,V,t) & \mapsto & (z,tV+\sum_{s\in 
    \fS_{sing}}Q_s(z,V);t,0) \end{array}$$ 

\noindent In the sequel we will give an idea of how one can show that there is a unique structure
of Lie
algebroid on $\cA$ such that $\Phi$ is its anchor map. The Lie algebroid $\cA$
is almost injective and so it is integrable, moreover we will see that
at a set level $\cG_X^t$ must be a groupoid which integrates it \cite{CF,D2001}. In
particular $\cG^t_X$ can be equipped with a unique smooth
structure such that it integrates the Lie algebroid $\cA$. 

\noindent Now we can state the following:

\begin{theorem} \label{PropLie algebroid} There exists a unique structure
  of Lie algebroid on the 
  smooth vector bundle $\cA=T\smx \times [0,1]$ over $\smx \times
  [0,1]$ with $\Phi$ as anchor.
\end{theorem}

\noindent To prove this theorem we will need several lemmas:

\begin{lemma}\label{lemQ} Let $s_0$ and $s_1$ be two strata such that
  $d(s_0)\leq d(s_1)$. \begin{enumerate} 
\item  For any tangent vector field $W$ on $\smx$,
  $Q_{s_1}(W)(\tau_{s_0})=0$. 
\item For any $(z,V)\in T\smx$, the following equality holds:
$$ Q_{s_1} \circ Q_{s_0}(z,V)= Q_{s_0} \circ
Q_{s_1}(z,V)=\tau_{s_0}(z) Q_{s_1}(z,V) \ .$$ 
\end{enumerate} 
\end{lemma}

\begin{proof} First notice that outside $O_{s_0}\cap O_{s_1}$
  either $Q_{s_1}$ hence $Q_{s_1}(W)$ or
  $\tau_{s_0}$ and $Q_{s_0}$ vanish 
thus the equalities in (1) and (2) are simply $0=0$. \\ 
(1) According to the compatibility conditions
\ref{compcond}  we have $\rho_{s_0}\circ
  \pi_{s_1}=\rho_{s_0}$ on
  $O_{s_0}\cap O_{s_1}$. Thus $\rho_{s_0}$ is constant on the fibers
  of $\pi_{s_1}$ and since $\tau_{s_0}=\tau \circ \rho_{s_0}$,
  $\tau_{s_0}$ is also constant on the fibers of 
  $\pi_{s_1}$. For any tangent vector field $W$, and any $z\in 
  \smO_{s_1}$ the vector  ${Q_{s_1}}(W)(z)$ is tangent to the fibers
  of $\pi_{s_1}$ thus $Q_{s_1}(V)(\tau_{s_0})=0$ on $O_{s_0}\cap O_{s_1}$.  \\
(2) The result follows from the first remark and the equality
\ref{tangent-filtration} of the part above.
\end{proof}

\noindent The next lemma ensures that $\Phi$ is almost injective, in particular it is
injective in restriction to $\smx\times ]0,1]$. A simple calculation
shows the following:

\begin{lemma}\label{pinjectif} For any $t\in ]0,1]$ the bundle map $\Phi_t$
  is bijective, moreover $$\Phi_t^{-1}(z)=\frac{1}{t}V- \sum_{s\in
    \fS_{sing}}
  \frac{1}{(t+t_s(z)) \cdot (t+t_{s}(z)-\tau_s(z))} Q_s(z,V) 
\ $$
where for any singular stratum $s$ the map $t_s$ is defined as follows:
  $$t_s: \smx \rightarrow \RR \ , \ t_s(z)=\sum_{s_0\leq s}^l
  \tau_{s_0}(z) \ .$$  
\end{lemma}

\noindent Thus in order to prove the theorem \ref{PropLie algebroid} it is enough to show
that locally the image of the map induced by $\Phi$ from the set of
smooth local sections of $\cA$ to the set of smooth local tangent vector fields
on $\smx\times [0,1]$ is stable under the Lie bracket.

\begin{proof}[Idea of the proof of Theorem \ref{PropLie algebroid}]
\smallskip \noindent First notice that outside the closure of $\cup_{s_i\in
  \fS_{sing}}\smO_{s_i}$ the image under ${\Phi}$ of local
tangent vector fields is clearly stable under
  Lie Bracket.\\ Thus using decomposition of the form
  \ref{tangent-graduation} described in the last part and standard
  arguments it remains to show that 
  if $s_a$ and $s_b$ are strata of depth respectively $a$ and $b$ with
  $s_a\leq s_b$, if $U$ is an open subset of $\smx$, as small as we
  want contained in $\cN_{s_a}\cap \cN_{s_b}$, and if
  $W^{\perp}, V^{\perp},\ V_a$ and $W_b$ are tangent vector fields on $U$, satisfying:
\begin{itemize} \item[] $V^{\perp}$ and $V_a$ can be porjected by $\pi_{s_b}$,
\item[] ${Q_s}(W^{\perp})={Q_s}(V^{\perp})=0 \ \mbox{for any } s\in 
  \fS$ ,
\item[] ${Q_s}(V_a)=\left\{ \begin{array}{ll}\tau_s V_a \mbox{ when }
      s\leq s_a \\ 0 \mbox{ elsewhere} \end{array} \right.$ and 
  ${Q_s}(W_a)=\left\{ \begin{array}{ll} \tau_s W_a \mbox{ when }
      s\leq s_b \\ 0 \mbox{ elsewhere} \end{array} \right.$,

\end{itemize}
then $[{\Phi}(W^{\perp}+W_b), {\Phi}(V^{\perp})]$ and 
$[{\Phi}(W_b), {\Phi}(V_a)]$ are 
in the image of ${\Phi}$. In other word, we have to show that the maps
$(z,t)\in \smX\times ]0,1] \mapsto (\Phi_t^{-1}([{\Phi}(W^{\perp}+W_b),
{\Phi}(V^{\perp})](z)),t) $ and $(z,t)\in \smX\times ]0,1] \mapsto
(\Phi_t^{-1}([{\Phi}(W_b), {\Phi}(V_a)](z)),t) $ can be extended into
smooth local section of $\cA$.
The result follows from our preceding lemmas and usual calculations.
\end{proof}

\noindent Now we can state:
\begin{theorem} The groupoid $\cG ^t_X$ can be equipped with a smooth
  structure such that its Lie algebroid is $\cA$ with $\Phi$ as
  anchor.
\end{theorem}

\begin{proof}
According to proposition \ref{PropLie algebroid} and lemma
\ref{pinjectif}, the Lie algebroid $\cA$ is almost injective. Thus
according to \cite{D2001} there is a unique
$s$-connected quasi-graphoid $\cG(\cA) \rightrightarrows \smx \times [0,1]$ which 
integrates $\cA$. Suppose for simplicity that for each stratum $s$,
$\smO_s$ is connected (which will ensure that $\cG^t_X\vert_{\smF_s\times
  [0,1]}$ is a $s$-connected quasi-graphoid). 

\noindent Moreover the map $\Phi$ satisfies:
\noindent \begin{enumerate}[\rm (i)] \item $\Phi$ induces an isomorphism from
  $\cA_{]0,1]}:=\cA\vert_{\smx\times ]0,1]}$ to $T\smx \times ]0,1]$, 
\item[\rm (ii)]
for any stratum $s$, the Lie algebroid $\cA$ restricted  over
$F_s^\circ\times [0,1]$ to a Lie
algebroid $\cA_s:=\cA\vert_{\smF_s\times [0,1]}$
which is isomorphic to the restriction of $\cA^t_{\pi_s\times \id}$
over $\smF_s\times [0,1]$. 
\end{enumerate}

\noindent Thus, again by using the uniqueness of $s$-connected quasi-graphoid integrating a
given almost injective Lie algebroid, we obtain:
\noindent \begin{enumerate}\item[\rm (i)] the restriction of the
  groupoid 
$\cG(\cA)$ over $\smx\times ]0,1]$ is isomorphic to $\smx \times \smx
\times ]0,1]\rightrightarrows \smx \times ]0,1]$, the pair groupoid
on $\smx$ parametrized by $]0,1]$, 
\item[\rm (ii)]  
for each stratum $s$ the restriction over 
$\smF_s\times [0,1]$ is equal to $\cG^t_X \vert_{\smF_s\times
  [0,1]}$. \end{enumerate}

\noindent Finally $\cG(\cA)=\cG^t_X$ and there is a unique smooth
structure on $\cG^t_X$ such that $\cA$ is its Lie algebroid.\\
If some $\smO_s$ is not connected, we replace in the construction of the
tangent space the groupoid $\pb{\pi_s}(Ts)\vert_{F_s}$ by its $s$-connected
component. Let $CT^{\fS}X$ and $C\cG^t_X$ be the corresponding
groupoids. The previous arguments apply and the groupoid $C\cG^t_X$ admits a unique smooth
structure such that $\cA$ is its Lie algebroid. One can then show that
there is a  unique smooth
structure on $\cG^t_X$ such that $C\cG^t_X$ is its $s$-connected
component. Precisely, according to \cite{D2001} there is a quasi-graphoid
$\cG\cI(\cA)\rightrightarrows \smX \times [0,1]$ which integrates $\cA$
and is maximal for the inclusion among quasi-graphoids which
integrate $\cA$. The groupoid $C\cG^t_X$ is then the $s$-connected
component of $\cG\cI(\cA)$. In particular it is open in
$\cG\cI(\cA)$. Let $X^r:= \smx \setminus
\underset{s\in \fS}\partial F_s$. The restriction of $\cG^t_X$ to $X^r \times [0,1]$ is a
quasi-graphoid which integrates the restriction of $\cA$ to $X^r \times [0,1]$ and is then clearly an open
sub-groupoid of $\cG\cI(\cA)$. Now we have $\cG^t_X=\{ \gamma\cdot \eta \
\vert \ \gamma \in C\cG^t_X,\ \eta \in \cG^t_X\vert_{X^r \times [0,1]},\ s(\gamma)=r(\eta)\}$
which is
open in $\cG\cI(\cA)$ and so $\cG^t_X$ inherits the required smooth structure.
\end{proof}

\noindent Thus
$T^{\fS}X$, which is the restriction of $\cG^t_X$ to the saturated set
$\smx \times \{0\}$, inherits from $\cG^t_X$ a smooth
structure which is equivalent to the one described in previous paragraphs.

\subsubsection{Standard projection from the tangent space onto the
  space}\label{proj} 
 
The space of orbits of $\smx /
  T^{\fS}X$ is equivalent to $X$ in the
sense that there is a canonical isomorphism  
$C_0(\smx /T^{\fS}X) \simeq C(X)$. 

\begin{definition} Let $r,s:T^\fS X\to X^\circ$ be the target and source maps of the
   $\fS$-tangent space of $X$. A continuous map $p:X\to X$ is a {\it standard projection } for
$T^{\fS}X$ on $X$ if:
\begin{enumerate}
\item  $\displaystyle p\circ r = p\circ s $.
\item  $p$ is homotopic to the identity map of $X$. 
\end{enumerate}
A standard projection $p$ for $T^{\fS}X$ on $X$ is {\sl surjective} if $p|_{\smx}:\smx\to X$ is onto. 
\end{definition}
This definition leads to the following:
\begin{lemma}\label{std-projection}
  \begin{enumerate}
  \item There exists a standard surjective projection for $T^{\fS}X$
    on $X$.  
  \item Two standard projections are homotopic and the homotopy can be
     done within the set of standard projections.
  \end{enumerate}
\end{lemma}
\begin{proof}
 1) If $X$ has depth $0$, $\smx=X$ and we just take $p=id$. Let us 
 consider $X$ with depth $k>0$. Choose a smooth non decreasing
 function 
$f:\RR_+\to\RR_+$ such that $f([0,1])=0$ and
$f|_{[2,+\infty[}=\mbox{Id}$. Recall that there exists for each singular
stratum $s$ an isomorphism \ref{mapingcylindre}: 
$$\Psi_s: \cN_s   \rightarrow 
    c_{\pi_s}L_s =L_s\times [0,+\infty[/\sim_s \ .$$
we define the map 
 $$
     p_s: \cN_s \longrightarrow \cN_s 
 $$
 by the formula:
 $$ 
    \Psi_s\circ p_s\circ \Psi_s^{-1}[x,t]= [x,f(t)].  
 $$ 
For each integer $i\in [0,k-1]$, we define a continuous map:
 $$ 
  p_i: X \longrightarrow X
 $$ 
by setting $p_i(z)=p_s(z)$ if $z$ belongs to $\cN_s$ for some singular
stratum of depth $i$ and $p_i(z)=z$ elsewhere. In particular,
$p_i|_{O_s}=\pi_s$ for every stratum $s$ of depth $i$. 
Finally we set:
 $$ 
  p = p_0\circ p_1\circ \cdots \circ p_{k-1}. 
 $$
This is the map we looked for. Indeed:

Let $\gamma\in T^\fS X$. There exists a unique stratum $s$ such
  that $\gamma\in \pb{\pi_s}(TS)$. If $s$ is regular, then 
$r(\gamma)=s(\gamma)$ so the result is trivial here. Let us assume
that $s$ is singular and let $i<k$ be its depth.   
By definition, $r(\gamma)$ and $s(\gamma)$ belong to $O_s$. For each
stratum $t \ge s $ of depth $j \ge i$, we have everywhere it makes sense: 
 $$
   \pi_t\circ p_t = \pi_t, \ \pi_s\circ\pi_t=\pi_s,\
   \rho_s\circ\pi_t=\rho_s
 $$
thus: 
 $$
  \rho_s\circ p_t=\rho_s\circ\pi_t\circ p_t=\rho_s\circ\pi_t=\rho_s
 $$
which proves that $p_j(O_s)=O_s$, and moreover:
 $$
  \pi_s\circ p_t= \pi_s\circ \pi_t\circ p_t=\pi_s\circ\pi_t=\pi_s
 $$
Recalling that
$p_i|_{O_s}=\pi_s|_{O_s}$, this last relation implies:
 $$ 
    p_i\circ \cdots \circ p_{k-1}|_{O_s} =\pi_s\circ p_{i+1}\circ
    \cdots \circ p_{k-1}|_{O_s} =\pi_s|_{O_s}
 $$
Since by definition we also have $\pi_s(r(\gamma))=\pi_s(s(\gamma))$,
we conclude that: 
 $$ 
   p(r(\gamma))=p_0\circ\cdots p_{i-1}\circ \pi_s(r(\gamma)) =
   p_0\circ\cdots p_{i-1}\circ \pi_s(s(\gamma))=  p(s(\gamma))
 $$
If in the definition of $p$,  we replace the function $f$ by
$t\id_{{\RR_{+}}} +(1-t)f$,  we get a homotopy between $p$ and $\id_{X}$. 

Finally,  $p$  has the required surjectivity property:
  $p_{k-1}(X^\circ)=X^\circ\bigcup_{d(s)=k-1}s$ and for all $j$ we
  have the equality 
$p_{j-1}(X^\circ\bigcup_{d(s)\ge j}s)=X^\circ\bigcup_{d(s)\ge j-1}s$. 

2) Let $q$ be a standard projection and $p$ be the standard
projection built in 1).  Let also $q_{t}$ be a homotopy
between $q$ and $\id_{X}$ and $p_{t}$ the homotopy built in 1)
between $p$ and $\id_{X}$. Observe that $q_{t}\circ p$ is a
standard projection,  providing a path of standard projections
between $q\circ p$ and $p$. Moreover,  by construction of $p_{t}$,
the inclusion $\im(p_{t}\circ r, p_{t}\circ s)\subset \im(r, s)$
holds for any $1\ge t>0$,  thus $q\circ p_{t}$
is a standard projection, providing a path of standard projections
between $q\circ p$ and $q$. Thus,  any standard projection $q$ is
homotopic to $p$ within the set of standard projections and the result
is proved. 
\end{proof}

\begin{remark}
Let $p$ be the surjective standard projection built in the proof of
the last proposition.  The map $p\circ r:T^{\fS}X\to X$ provides
$T^{\fS}X$ with a structure of continuous field of
groupoids. Following the arguments of (\cite{DL2},  remark 5), it can be shown
that each fiber of this field is amenable,  thus $T^{\fS}X$ is
amenable and $C^{*}(T^{\fS}X)=C^{*}_{r}(T^{\fS}X)$ is nuclear.  The
same holds for $\cG^{t}_{X}$ and all other deformation groupoids used below. 
\end{remark}


\section{Poincar{\'e} duality for stratified pseudo-manifolds}

Let $X$ be a compact stratified pseudomanifold of depth $k\ge 0$. 

The tangent groupoid
$\cG^t_X$ is a deformation groupoid,  thus it provides us with a
$K$-homology class, called a {\it pre-Dirac} element:
\begin{equation}
  \label{prediracX}
  \delta_X= [e_0]^{-1}\otimes [e_1]\in KK(C^*(T^\fS X),\CC).
\end{equation}
Here $e_0: C^*(\cG^t_X)\to C^*(T^\fS X)$ and 
$e_1: C^*(\cG^t_X)\to\cK(L^2(\smx))$ are the usual evaluation
homomorphisms. Now we need:
\begin{lemma}  1) Let $p:\smx \rightarrow X$ be a surjective standard projection for $T^{\fS}X$. The formula:
 $$
   \forall a\in C^*(T^\fS X),\ f\in C(X),\  \gamma\in T^\fS
   X,\quad  (a\cdot f)(\gamma)=f(p\circ r(\gamma)).a(\gamma)
 $$
defines a $C(X)$-algebra structure on $C^*(T^\fS X)$.\\ 
 2) For any standard projection $p$ for $T^{\fS}X$, the formula: 
 $$
   \forall a\in C^*(T^\fS X),\ f\in C(X),\  \gamma\in T^\fS
   X,\quad  \Psi_X(a\cdot f )(\gamma)=f(p\circ r(\gamma)).a(\gamma)
 $$
defines a homomorphism $\Psi_X:   C^*(T^\fS X) \otimes C(X)  \to C^*(T^\fS X)$ 
whose class $[\Psi_X]\in KK(C^*(T^\fS X)\otimes C(X),C^*(T^\fS X) )$
does not depend on the choice of $p$.
\end{lemma}
The last assertion uses Lemma \ref{std-projection}. Note that if
$k=0$, $X$ is smooth and we can choose $p=\id$,  thus: 
\begin{equation}
  \label{psi-smooth}
     \Psi_X(a\otimes b)(V) = b(x).a(x,V)
\end{equation}
for all $V\in T_xX$, $x\in X$,  $a\in C(X)$ and $b\in C^{*}(TX)$. 
 
From now on, we choose a surjective standard projection and denote by 
$\Psi_{X}: C(X) \otimes C^*(T^\fS X)  \to C^*(T^\fS X)$ the homomorphism defined in the previous lemma. We set:
\begin{equation}
  \label{diracX}
   D_X = \Psi^*_X(\delta_X)=[\Psi_X]\otimes \delta_X\in KK(C^*(T^\fS X)\otimes C(X),\CC).
\end{equation}
This section is devoted to the proof of the main theorem:
                                                                                                      
\begin{theorem}\label{main-theorem} Let $X$ be a compact stratified pseudomanifold. The
  $K$-homology class $D_X$ is a Dirac element, that is, it provides a
  Poincar\'e duality between the algebras $C^*(T^\fS X)$ and $C(X)$. 
\end{theorem}

We need some notations. If $W$ is an open set of the stratified pseudomanifold $X$ and
$\overline{W}$ its closure, we set:
$$T^\fS W= T^\fS X|_{W^{\circ}}\quad ; \quad  T^\fS \overline{W}=T^\fS
X|_{\overline{W}^\circ} \quad  \mbox{ and
} \quad
\cG^t_W=\cG^t_X\vert_{W^{\circ}\times [0,1]}. $$
The groupoid $\cG^t_W$ is a deformation groupoid which defines the
$K$-homology class $\delta_W \in K^0(C^*(T^\fS W))$.
We define the homomorphisms  induced by
$\Psi_X$: $$\hat{\Psi}_W:C^*(T^\fS
\overline{W})\otimes C_0(W) \rightarrow C^*(T^\fS
{W}) \quad \mbox{ and
} \quad \hat{\Psi}_{\overline{W}}:C^*(T^\fS
{W})\otimes C(\overline{W}) \rightarrow C^*(T^\fS
{W})\, $$ and we set ${\Psi}_W=i_W\circ \hat{\Psi}_W$ and
${\Psi}_{\overline{W}}=i_W\circ \hat{\Psi}_{\overline{W}}$ where
$i_w : C^*(T^\fS
{W}) \rightarrow C^*(T^\fS X)$ is the natural homomorphism. Finally we
    let : $$D_W=(\hat{\Psi}_W)^*(\delta_W)=(\Psi_W)^*(\delta_X)
    \in KK(C^*(T^\fS \overline{W})\otimes C_0(W),\CC)$$
and 
$$D_{\overline{W}}=(\hat{\Psi}_{\overline{W}})^*(\delta_W)=(\Psi_{\overline{W}})^*(\delta_X)
    \in KK(C^*(T^\fS W)\otimes C(\overline{W}),\CC)\ .$$
In the sequel,  we will be interested in the disjoint open sets:
\begin{equation}
  \label{eq:O-pm}
   O_{-}=\bigcup_{s\in \fS_0}\{z\in \cN_{s} \ \vert \ \rho_s(z) <2\}
   \hbox{ and }
       O_{+}= X\setminus \overline{O_{-}},
\end{equation}
as well as in the intersection of their closures:
\begin{equation}
  \label{eq:L}
  L=\overline{O_{+}}\cap \overline{O_{-}}=\bigcup_{s\in \fS_0} \{z\in X\ \vert\ \rho_s(z)=2\}.
\end{equation}
We recall from Paragraph (\ref{unfolding}) that $\fS_0$ denotes the set of
minimal strata. 

\begin{proof}[Proof of Theorem \ref{main-theorem}]
\smallskip \noindent It will be proved by induction on the depth
of the stratification and the unfolding process will be used to reduce the depth. 

\smallskip \noindent If $\mbox{depth}(X)=0$ the content of the theorem is well known, and
that $D_X$ is a Dirac element is a consequence of \cite{DL2}. 
Let $k\ge 0$, assume that the theorem \ref{main-theorem} holds for all compact
stratified pseudomanifolds with depth $\le k$ and let $X$ be a compact
stratified pseudomanifold of depth 
$k+1$. The proof of the induction is divided in two parts. 

{\bf First part of the proof.}
We consider two natural ``restrictions'' of $D_{X}$, namely 
$D_{O_+}\in K^{0}(C^*(T^{\fS}\overline{O_{+}})\otimes C_0(O_+))$ and 
$D_{\overline{O_-}}\in K^{0}(C^*(T^{\fS}O_{-})\otimes C(\overline{O_{-}}))$. Then, we reduce
the proof of the theorem to checking that $D_{O_{+}}$ is a Dirac
element.  

Let $O_0$ be the open set of $X$ obtained by replacing the
condition $\rho_{s}<2$ by $\rho_{s}<1$ in the definition of $O_{-}$
in (\ref{eq:O-pm}). The $C^{*}$-algebra $C^{*}(T^{\fS}O_0)$ is a
closed two-sided ideal in $C^{*}(T^{\fS}O_{-})$ and the quotient 
 $$
   C^{*}(T^{\fS}O_{-})/C^{*}(T^{\fS}O_0)\simeq C_{0}([1,2[)\otimes C^{*}(T^{\fS}L\times\RR)
 $$
is contractible in $K$-theory. It follows that the inclusion 
$C^{*}(T^{\fS}O_0)\subset C^{*}(T^{\fS}O_{-})$ is a
$KK$-equivalence which sends $\delta_{O_{0}}$ to $\delta_{O_-}$.
 This is obvious once we consider the corresponding tangent groupoids 
$\cG^{t}_{O_{0}},\ \cG^{t}_{O_-}$. As already noted there is a natural Morita
equivalence between the groupoid
$T^{\fS}O_0=\bigcup_{s\in\fS_0}\pb{\pi_s}(Ts)|_{O_{s}}$ and the
tangent space $TS= \bigcup_{s\in\fS_0}Ts$ of the closed smooth
manifold $S=\cup_{s\in\fS_{0}}s$. Under this Morita equivalence,
$\delta_{O_{0}}$ corresponds to $\delta_{S}$: this follows from the
extension of the previous Morita equivalence to the tangent groupoids 
$\cG^{t}_{O_{0}}$ and $\cG^{t}_{S}$. Moreover the control data
provide a homotopy equivalence betwen $\overline{O_{-}}$ and $S$ and we finally
get a $KK$-equivalence between $C^{*}(T^{\fS}O_{-})\otimes C(\overline{O_{-}})$
and $C^{*}(T^{\fS}S)\otimes C(S)$ under which the class
$[\Psi_{\overline{O_{-}}}]$ coincides with the class $[\Psi_{S}]$. We have proved: 
\begin{lemma}\label{mthm-prop-2}
  There is a $KK$-equivalence  between $C^{*}(T^{\fS}O_{-})\otimes
  C(\overline{O_{-}})$ and $C^{*}(T^{\fS}S)\otimes C(S)$ under which the Dirac element $D_S$ corresponds to $D_{\overline{O_{-}}}$. 
  In particular, $D_{\overline{O_{-}}}$ is a Dirac element.
\end{lemma}
We now apply Lemma \ref{2-ou-of-3} to the nuclear 
$C(X)$-algebra $C^*(T^\fS X)$, the disjoint open subsets $O_{-}$ and
$O_+$ and the $K$-homology class $\delta_{X}$. 
Since by Lemma \ref{mthm-prop-2} $D_{\overline{O_{-}}}$ is a Dirac element,  we
immediately get:
\begin{center}
 $D_X$ is a Dirac element if and only if $D_{O_+}$ is. 
\end{center}

{\bf Second part of the proof.}
We check that $D_{O_{+}}$ is a Dirac element.  
Let us go back to the compact pseudomanifold of depth $k$ coming from the unfolding
process: $2X$. Modifying slightly the definition of Paragraph \ref{unfolding},  we set: 
 $$
   2X=\overline{O_+} \underset{L}{\cup} L\times [-2,+2] \underset{L}{\cup} \overline{O_+}
 $$
We consider this time the disjoint open subsets $U=L\times ]-2,+2[$
and $V=2X\setminus L\times [-2,+2]=O_{+}\sqcup O_{+}$ of the
pseudomanifold $2X$. 
Let us introduce as before the homomorphisms induced by $\Psi_{2X}$:
\begin{equation}
  \label{psiLR}
   \Psi_{\bU}: C^*(T^\fS U)\otimes
   C(\bU) \longrightarrow C^*(T^\fS 2X)
\end{equation}
and 
\begin{equation}
  \label{psipm}
   \Psi_{V}: C^*(T^\fS\bV)\otimes  C_{0}(V) \longrightarrow C^*(T^\fS 2X)
\end{equation}
where $T^\fS U=T^\fS 2X|_{U^{\circ}}$ and $T^\fS\bV=T^\fS 2X|_{\bV^{\circ}}$. Note that under the natural identification 
$C^*(T^\fS\bV)\otimes  C_{0}(V)\simeq
M_{2}(C^{*}(T^{\fS}\overline{O_{+}})\otimes C_{0}(O_{+}))$,  the homomorphism
$\Psi_{V}$ has the  following diagonal form: $\Psi_{V}=\mathrm{diag}(\Psi_{O_{+}},\Psi_{O_{+}})$. 

We shall consider three $K$-homology classes: 
 $$
   D_{2X}=\Psi^{*}_{2X}(\delta_{2X}), \quad
   D_{\bU}=\Psi^{*}_{\bU}(\delta_{2X}), \quad
   D_{V}=\Psi^{*}_{V}(\delta_{2X})\ .
 $$
\noindent Since $2X$ is a compact stratified pseudomanifold of depth $k$, we know
by induction hypothesis that $D_{2X}$ is a
Dirac element. 

The space $L$ with the stratification induced by $X$ is also a compact 
stratified pseudomanifold of depth $k$. So it has a Dirac element
$D_{L}$ defined as before. Observe that $\Psi_{\bU}$ has
range in the ideal $C^*(T^\fS U)$ of $C^*(T^\fS 2X)$.  We note
$\hat{\Psi}_{\bU}$ the induced homomorphism,  $i_U:C^*(T^\fS U)\to C^*(T^\fS
2X)$ the inclusion and $\delta_{U}$
the $KK$-element associated with the deformation groupoid 
$\cG^{t}_{U}:=\cG^{t}_{2X}|_{U^{\circ}}$. We have
$\delta_{U}=(i_U)^{*}(\delta_{2X})$,  hence $D_{\bU} =
(\hat{\Psi})^{*}_{\bU}(\delta_{U})$. On the other hand, let $\delta$
be the $KK$-element associated with the deformation groupoid
$\cG^{t}_{]-2, 2[}$.  It is clear that $\delta$ is a
generator of $K^{0}(C^{*}(T]-2, 2[))\simeq\ZZ$ and its pull-back
$\Delta$ under the homotopy equivalence $C([-2,2])\to\CC$ is a Dirac
element. Now, under the groupoid
isomorphism $T^{\fS}U\simeq  T^{\fS}L \times T]-2, 2[$,  the element $\delta_{U}$
corresponds to $\delta_{L}\underset{\CC}{\otimes}\delta$ and 
$D_{\bU}$ to $D_{L}\underset{\CC}{\otimes}\Delta$. It follows that
$D_{\bU}$ is a Dirac element.

Since $D_{2X}$ and $D_{\bU}$ are Dirac elements, we get from Lemma
\ref{2-ou-of-3} applied to the nuclear
$C(2X)$-algebra $C^{*}(T^{\fS}2X)$,  to the open sets $U, V$ and to the $K$-homology
class $\delta_{2X}$,  that $D_{V}$ is a Dirac element. Since
$\Psi_{V}^{*}$ has diagonal form, 
we have:
\begin{equation}
  \label{double-DOb}
    D_{V}=D_{O_+}\oplus D_{O_+} \in K^0(C^{*}(T^{\fS}\overline{O_{+}})\otimes 
    C_{0}(O_{+}))^{\oplus^{2}}\subset K^0(C^*(T^\fS\bV)\otimes C_0(V)).  
\end{equation}
It is clear from this formula that $D_{V}$ is a Dirac element if and
only if $D_{O_{+}}$ is,  so we have proved that $D_{O_{+}}$ is a
Dirac element,  which ends the proof of the theorem.  
\end{proof}

The following remark collect some technical facts which
were in the main body of the proof of Theorem  \ref{main-theorem} before we
took into account the Referee's suggestions.
\begin{remark}
Let us replace $\rho_{s}(z)<2$ by $\rho_{s}(z)<1$ in the
   definition of $O_{-}$ in (\ref{eq:O-pm}) and modify accordingly
   the subsequent sets in (\ref{eq:O-pm}). 
   Let $\partial_{\pm}\in KK_{1}(C^{*}(T^{\fS}L\times\RR),C^{*}(T^{\fS}O_{\pm}))$ 
      be the $KK$-elements associated with the exact sequences of $C^{*}$-algebras:
\begin{equation}
  \label{exact-pm}
  0\longrightarrow C^{*}(T^{\fS}O_{\pm})\longrightarrow
  C^{*}(T^{\fS}\overline{O_{\pm}})\longrightarrow
  C^{*}(T^{\fS}L\times\RR)\longrightarrow 0. 
\end{equation}
We can apply Lemma \ref{referee-1} to $A=C^{*}(T^{\fS}X)$,
$J_{1}=C^{*}(T^{\fS}O_{+})$ and $J_{2}=C^{*}(T^{\fS}O_{-})$. This
gives: 
 \begin{equation}
  \label{rest-O-pm}
  \partial_{+}\otimes \delta_{O_{+}} = - \partial_{-}\otimes
  \delta_{O_{-}}\in K^{1}(C^{*}(T^{\fS}L\times\RR)).
 \end{equation}
Moreover, one can show that this element is, modulo sign and Bott periodicity 
$K^{1}(C^{*}(T^{\fS}L\times\RR))\simeq K^0(C^{*}(T^{\fS}L))$, the Dirac element $D_{L}$ associated
with $L$.  The idea to prove this is to build a smooth groupoid: 
 $$ 
  \hat{\cG}^t_{\overline{O_{-}}}:= \cG^t_{O_{-}}\sqcup (\cG^t_L\times \RR) \rightrightarrows \overline{O_-}\times [0,1].
 $$
such that the following (smooth) isomorphisms hold: 
\begin{itemize}
 \item $\displaystyle \hat{\cG}^t_{\overline{O_{-}}}|_{\overline{O_-}\times\{0\}}\simeq T^\fS \overline{O_{-}}$,
 \item $\displaystyle \hat{\cG}^t_{\overline{O_{-}}}|_{\overline{O_-}\times\{1\}}\simeq 
  (L^\circ\times L^\circ)\times (\RR\rtimes_\phi\RR|_{]0,1]}) $,
\end{itemize}
where  $ \RR\rtimes_\phi\RR\rightrightarrows \RR$ is the groupoid of the action of $\RR$ onto itself by
the complete flow of the vector field $\tau(h)\partial_h$ and $\tau$
is the gluing function used in Paragraph
\ref{subsection:smooth-structure} (this is exactly the tangent space
of $[1, +\infty[$ with $\{1\}$ as a conical point). Since
$\hat{\cG}^t_{\overline{O_{-}}}|_{\overline{O_-}\times\{1\}}$ has vanishing  
$K$-theory, hence the $KK$-element $\alpha$ associated with the exact sequence:
\begin{equation}
  \label{exact-seq-bert}
  0\longrightarrow C^*(O_{-}^\circ\times O_{-}^\circ) \longrightarrow
  C^*(\hat{\cG}^t_{\overline{O_{-}}}|_{\overline{O_-}\times\{1\}})\longrightarrow C^*(L^\circ\times
  L^\circ\times\RR)\longrightarrow 0
\end{equation}
is invertible in $KK$-theory,  thus corresponds to Bott periodicity
modulo a sign and the  Morita equivalences between 
$C^{*}(O_{-}^\circ\times O_{-}^\circ), C^*(L^\circ\times L^\circ) $
and $\CC$. Finally, we consider the commutative diagram:  
\begin{equation}
  \label{proof-lem}
  \xymatrix{
  0 \ar[r] & C^*(T^\fS O_{-}) \ar[r] & C^*(T^\fS \overline{O_{-}})\ar[r] & C^*(T^\fS
  L\times\RR) \ar[r] & 0 \\
   0 \ar[r] & C^*(\cG^t_{O_{-}}) \ar[u]_{e_0^{O_{-}}}\ar[d]^{e_1^{O_{-}}}\ar[r] &
   C^*(\hat{\cG}^t_{\overline{O_{-}}})\ar[u]_{e_0^{\overline{O_{-}}}}\ar[d]^{e_1^{\overline{O_{-}}}}\ar[r] & 
   C^*(\cG^t_{L}\times\RR)\ar[u]_{e_0^L\otimes 1} \ar[d]^{e_1^L\otimes 1} \ar[r] & 0 \\
  0 \ar[r] & C^*(O_{-}^\circ\times O_{-}^\circ) \ar[r] &
  C^*(\hat{\cG}^t_{\overline{O_{-}}}|_{\overline{O_-}\times\{1\}})\ar[r] & C^*(L^\circ\times
  L^\circ\times\RR) \ar[r] & 0 }
\end{equation}
It gives by functoriality: 
$\displaystyle\partial_{-}\otimes\delta_{O_{-}}=\delta_L\underset{\CC}{\otimes}
\alpha$ which proves the claim. 

\end{remark}

\subsubsection{Stratified pseudomanifold with boundary.} 
As a
byproduct of the proof of Theorem \ref{main-theorem}, we have proved that Poincar\'e duality 
also holds for compact stratified pseudomanifolds with
boundary. Precisely a stratified pseudomanifold with boundary is
$(X_b,L,\fS_b,N_b)$ where:
\begin{enumerate} \item $X_b$ is a compact 
  separable metrizable space and
  $L$ is a compact subspace of $X_b$.
\item $\fS_b=\{s_i\}$ is a finite partition of $X_b$ into locally closed
  subset of $X_b$, which are smooth manifolds possibly with 
  boundary. Moreover for each $s_i$ we have $$s_i\cap L =\partial s_i
  \ .$$ 
\item $N_b=\{\cN_{s},\pi_{s},\rho_{s}\}_{s\in \fS_b}$, where
  $\cN_s$ is an open neighborhood of $s$ in 
  $X$, $\pi_s:\cN_s \rightarrow s$ is a continuous retraction and
  $\rho_s:\cN_s \rightarrow [0,+\infty[$ is a continuous map such that
  $s=\rho_s^{-1}(0)$.
\item The {\it double}: $$X=X_b\underset{L}{\cup} X_b$$ obtained by
  gluing two copies of $X_b$ along $L$ together with the partition
  $\fS:=\{s_i\ \vert \partial s_i=\emptyset \} \cup
\{s_i\underset{\partial s_i} \cup s_i\} \cup \{s_i\ \vert \partial
s_i=\emptyset \}$ and the set of 
  control data $N=\{ \tilde{\cN}_{s},\tilde{\pi}_{s},\tilde{\rho}_{s}
  \}_{s\in \fS}$  where $$\cN_s=\cN_{s_i},\ \pi_s=\pi_{s_i}, \
  \rho_s=\rho_{s_i} \mbox{ if } s=s_i \mbox{ with } \partial
  s_i=\emptyset $$ and $$\cN_s=\cN_{s_i}\underset{\cN_{s_i}\cap L}{\cup}
  \cN_{s_i},\ \pi_s\vert _{\cN_{s_i}\setminus L}=\pi_{s_i}, \
  \rho_s\vert _{\cN_{s_i}\setminus L} =\rho_{s_i} \mbox{ elsewhere }$$
  is a stratified pseudomanifold.
\end{enumerate}

\noindent We let $O_b:= X_b\setminus L$. According to the previous work, one can
define the tangent spaces: $$T^{\fS}X_b:= T^{\fS} X \vert_{X_b}
\mbox{ and } T^{\fS} O_b:= T^{\fS} X\vert_{O_b} $$
We deduce the following:
\begin{theorem} The $C^*$-algebras $C^*(T^{\fS}X_b)$ and $C_0(O_b)$ are
  Poincar{\'e} Dual as well as the $C^*$-algebras $C^*(T^{\fS}O_b)$ and
  $C(X_b)$.
\end{theorem}


\begin{thebibliography}{10}

\bibitem{AR2000}
C.~Anantharaman-Delaroche and J.~Renault.
\newblock {\em Amenable groupoids}, volume~36 of {\em Monographies de
  L'Enseignement Math\'ematique [Monographs of L'Enseignement Math\'ematique]}.
\newblock L'Enseignement Math\'ematique, Geneva, 2000.
\newblock With a foreword by Georges Skandalis and Appendix B by E. Germain.

\bibitem{AS1}
M.~Atiyah and I.~Singer.
\newblock The index of elliptic operators {I}.
\newblock {\em Annals of Math.}, 87:484--530, 1968.

\bibitem{BHS}
J.-P. Brasselet, G.~Hector, and Saralegi.
\newblock Th\'eor\`eme de de {R}ham pour les vari\'et\'es stratifi\'ees.
\newblock {\em Ann. Global Anal. Geom.}, 9(3):211--243, 1991.

\bibitem{BHS1992}
J.P. Brasselet, G.~Hector, and M.~Saralegi.
\newblock {$\mathcal{L}^2$}-cohomologie des espaces stratifi\'es.
\newblock {\em Manuscripta Math.}, 76(1):21--32, 1992.

\bibitem{CW}
A.~{C}annas {da}~Silva and A.~Weinstein.
\newblock {\em Geometric Models for Noncommutative Algebras}.
\newblock Berkeley Math. Lecture Notes series, 1999.

\bibitem{Co3}
A.~Connes.
\newblock A survey of foliations and operators algebras.
\newblock In providence AMS, editor, {\em Operator algebras and applications,
  Part 1}, volume~38 of {\em Proc. Sympos. Pure Math.}, pages 521--628, 1982.

\bibitem{Co0}
A.~Connes.
\newblock {\em Noncommutative {G}eometry}.
\newblock Academic Press, 1994.

\bibitem{AC-GS1984}
A.~Connes and G.~Skandalis.
\newblock The longitudinal index theorem for foliations.
\newblock {\em Publ. R.I.M.S. Kyoto Univ.}, 20:1139--1183, 1984.

\bibitem{CF}
M.~Crainic and R.L. Fernandes.
\newblock Integrability of {L}ie brackets.
\newblock {\em An. of Math.}, 157:575--620, 2003.

\bibitem{JC-GS1986}
J.~Cuntz and G.~Skandalis.
\newblock Mapping cones and exact sequences in {$KK$}-theory.
\newblock {\em J. Operator Theory}, 15(1):163--180, 1986.

\bibitem{D2001}
C.~Debord.
\newblock Holonomy groupoids for singular foliations.
\newblock {\em J. of Diff. Geom.}, 58:467--500, 2001.

\bibitem{DL2007}
C.~Debord and J.-M. Lescure.
\newblock Index theory and groupoids.
\newblock Notes of the lectures given at the summer school {\it Geometric and
  Topological Methods for Quantum Field Theory} at Villa de Leyva in July 2007.

\bibitem{DL2}
C.~Debord and J.-M. Lescure.
\newblock $k$-duality for pseudomanifolds with isolated singularities.
\newblock {\em J. Functional Analysis}, 219(1):109--133, 2005.

\bibitem{DLN2006}
C.~Debord, J.M. Lescure, and V.~Nistor.
\newblock Groupoids and an index theorem for conical pseudo-manifolds.
\newblock Preprint arxiv:0609438. To appear in J.
  f{\"u}r die reine und angewandte Mathematik, 2006.

\bibitem{EM2007}
H.~Emerson and R.~Meyer.
\newblock Dualities in equivariant {K}asparov theory.
\newblock Preprint arxiv:0711.0025.

\bibitem{EM2005}
H.~Emerson and R.~Meyer.
\newblock Euler characteristics and {G}ysin sequences for group actions on
  boundaries.
\newblock {\em Math. Ann.}, 334(4):853--904, 2006.

\bibitem{GoMa}
M.~Goresky and R.~MacPherson.
\newblock {I}ntersection homology theory.
\newblock {\em Topology}, 19:135--162, 1980.

\bibitem{HS1}
M.~Hilsum and G.~Skandalis.
\newblock Morphismes ${K}$-orient\'{e}s d'espaces de feuilles et
  fonctorialit\'{e} en th\'{e}orie de {K}asparov.
\newblock {\em Ann. Sci. Ecole Norm. Sup.}, 20{\ }(4):325--390, 1987.

\bibitem{HW}
Bruce Hughes and Shmuel Weinberger.
\newblock Surgery and stratified spaces.
\newblock In {\em Surveys on surgery theory, Vol. 2}, volume 149 of {\em Ann.
  of Math. Stud.}, pages 319--352. Princeton Univ. Press, Princeton, NJ, 2001.

\bibitem{Ka1}
G.G. Kasparov.
\newblock The operator {K}-functor and extensions of ${C}^*$-algebras.
\newblock {\em Izv. Akad. Nauk SSSR, Ser. Math.}, 44:571--636, 1980.

\bibitem{Ka2}
G.G. Kasparov.
\newblock Equivariant {KK}-theory and the {N}ovikov conjecture.
\newblock {\em Invent. math.}, 91:147--201, 1988.

\bibitem{JM_L2006}
Jean-Marie Lescure.
\newblock Elliptic symbols, elliptic operators and {P}oincar\'e duality on
  conical pseudomanifolds.
\newblock Preprint. arXiv:math.OA/0609328. To appear in J. of
  K-Theory.

\bibitem{Mc}
K.~Mackenzie.
\newblock {\em Lie groupoids and Lie algebroids in differential geometry},
  volume 124 of {\em London Mathematical Society Lecture Note}.
\newblock Cambridge university press, 1987.

\bibitem{Mat}
John~N. Mather.
\newblock Stratifications and mappings.
\newblock In {\em Dynamical systems (Proc. Sympos., Univ. Bahia, Salvador,
  1971)}, pages 195--232. Academic Press, New York, 1973.

\bibitem{M1990}
R.~Melrose.
\newblock Pseudodifferential operators, corners and singular limits.
\newblock In {\em Proceedings of the International Congress of
  Mathematicians,Vol.\ I, I (Kyoto, 1990)}, pages 217--234, Tokyo, 1991. Math.
  Soc. Japan.

\bibitem{MR2006}
R.~Melrose and F.~Rochon.
\newblock Index in {K}-theory for families of fibred cusp operators.
\newblock Preprint. arXiv:math.DG/0507590v2.

\bibitem{Mo1}
B.~Monthubert.
\newblock {\em Groupo\"{\i}des et calcul pseudo-diff\'{e}rentiel sur les
  vari\'{e}t\'{e}s \`{a} coins}.
\newblock PhD thesis, Universit\'{e} Paris {VII}-Denis Diderot, 1998.

\bibitem{MN2005}
B.~Monthubert and V.~Nistor.
\newblock A topological index theorem for manifolds with corners.
\newblock Preprint. arXiv:math.KT/0507601.

\bibitem{MRW}
P.~Muhly, J.~Renault, and D.~Williams.
\newblock Equivalence and isomorphism for groupoid {$C\sp \ast$}-algebras.
\newblock {\em J. Operator Theory}, 17(1):3--22, 1987.

\bibitem{NSS2}
V.~E. Naza{\u\i}kinski{\u\i}, A.~Yu. Savin, and B.~Yu. Sternin.
\newblock On the homotopy classification of elliptic operators on stratified
  manifolds.
\newblock {\em Dokl. Akad. Nauk}, 408(5):591--595, 2006.

\bibitem{Re}
J.~Renault.
\newblock {\em A groupoid approach to $C^*$-algebras}, volume 793 of {\em
  Lecture Notes in Math.}
\newblock Springer-Verlag, 1980.

\bibitem{Sav1}
A.~Savin.
\newblock Elliptic operators on manifolds with singularities and
  {$K$}-homology.
\newblock {\em $K$-Theory}, 34(1):71--98, 2005.

\bibitem{Ver}
A.~Verona.
\newblock {\em Stratified mappings---structure and triangulability}, volume
  1102 of {\em Lecture Notes in Mathematics}.
\newblock Springer-Verlag, Berlin, 1984.

\bibitem{W1965}
H.~Whitney.
\newblock Local properties of analytic varieties.
\newblock In {\em Differential and {C}ombinatorial {T}opology ({A} {S}ymposium
  in {H}onor of {M}arston {M}orse)}, pages 205--244. Princeton Univ. Press,
  Princeton, N. J., 1965.

\end{thebibliography}
 

%
%

\bibliographystyle{gtart}

\end{document}